\documentclass[12pt,a4paper,reqno]{amsart}
\usepackage[utf8]{inputenc}
\usepackage{amsmath}
\usepackage{amsfonts}
\usepackage{amssymb}
\usepackage{amsmath,bbm,amssymb,amsxtra}
\usepackage{mathrsfs}
\usepackage{enumerate}
\usepackage{caption}
\allowdisplaybreaks
\usepackage{epsfig}
\usepackage{ae}
    \def\Xint#1{\mathchoice
    {\XXint\displaystyle\textstyle{#1}}%
    {\XXint\textstyle\scriptstyle{#1}}%
    {\XXint\scriptstyle\scriptscriptstyle{#1}}%
    {\XXint\scriptscriptstyle\scriptscriptstyle{#1}}%
    \!\int}

    \def\XXint#1#2#3{{\setbox0=\hbox{$#1{#2#3}{\int}$ }
    \vcenter{\hbox{$#2#3$ }}\kern-.6\wd0}}
    \def\dashint{\Xint-}

\theoremstyle{definition}

\newtheorem{lemma}{Lemma}[section]

\newtheorem{proposition}[lemma]{Proposition}

\newtheorem{theorem}[lemma]{Theorem}

\newtheorem{corollary}[lemma]{Corollary}

\newtheorem{remark}[lemma]{Remark}

\newtheorem{definition}[lemma]{Definition}
\newtheorem*{notation}{Notation}

\newcommand{\prop}[1]{\begin{proposition}\label{#1} \sl }
\newcommand{\eprop}{\end{proposition}}
\newcommand{\thm}[1]{\begin{theorem}\label{#1} \sl }
\newcommand{\ethm}{\end{theorem}}
\newcommand{\cor}[1]{\begin{corollary}\label{#1} \sl }
\newcommand{\ecor}{\end{corollary}}

\newcommand{\lem}[1]{\begin{lemma}\label{#1} \sl }
\newcommand{\elem}{\end{lemma}}
\newcommand{\defin}[1]{\begin{definition}\label{#1} \sl }
\newcommand{\edefin}{\end{definition}}
\newcommand{\beqno}{\begin{eqnarray*}}
\newcommand{\eeqno}{\end{eqnarray*}}
\newcommand{\beqla}[1] {\begin {eqnarray}\label{#1}}
\def\eeq {\end {eqnarray}}
\newcommand{\beq}{\begin {eqnarray}}
\newcommand{\real}{{\mathbb R}}

\newcommand{\integer}{{\mathbb Z}}
\newcommand{\R}{{\mathbb R}}
\newcommand{\nanu}{{\mathbb N}}

\newcommand{\complex}{{\mathbb C}}

\newcommand{\F}{{\mathcal F}}

\renewcommand{\:}{\colon}

\newcommand{\res}{\mathcal R}
\newcommand{\ext}{\mathcal E}

\newcommand{\Q}{{\mathbf Q}}

\newcommand{\dist}{{\rm dist}\,}
\newcommand{\sdist}{{\rm dist}^*}
\newcommand{\supp}{{\rm supp}\,}
\newcommand{\Lip}{{\rm Lip}\,}

\newcommand{\hlmax}{\mathcal{M}}

\newcommand{\locint}{L^1_{\rm{loc}}}
\newcommand{\B}{\mathcal{B}}
\newcommand{\bp}{\B^s_{p,q}(\real^{d},\mu)}
\newcommand{\bpd}{\B^s_{p,q}(\real^d)}
\newcommand{\bpa}{\B^s_{p,q}(\real^{d+1}_+,\mu_\alpha)}
\newcommand{\fp}{\F^s_{p,q}(\real^{d},\mu)}
\newcommand{\fpd}{\F^s_{p,q}(\real^d)}
\newcommand{\fpa}{\F^s_{p,q}(\real^{d+1}_+,\mu_\alpha)}
\newcommand{\bpt}{\B^{s - (\alpha+1)/p}_{p,q}(\real^d)}
\newcommand{\bpad}{\B^s_{p,p}(\real^{d+1}_+,\mu_\alpha)}

\newcommand{\bptd}{\B^{s - (\alpha+1)/p}_{p,p}(\real^d)}
\newcommand{\dyadic}{\mathscr{Q}}
\newcommand{\W}{\mathscr{W}}
\newcommand{\A}{\mathscr{A}}
\newcommand{\U}{\mathscr{U}}
\newcommand{\N}{\mathscr{N}}
\newcommand{\T}{\mathscr{T}}

\newcommand{\bps}{\B^s_{p,q}(\real^{d+n},\mu)}
\newcommand{\fps}{\F^s_{p,q}(\real^{d+n},\mu)}
\newcommand{\bpas}{\B^s_{p,q}(\real^{d+n},\mu_\alpha)}
\newcommand{\fpas}{\F^s_{p,q}(\real^{d+n},\mu_\alpha)}
\newcommand{\bpts}{\B^{s - (\alpha+n)/p}_{p,q}(\real^d)}
\newcommand{\bptds}{\B^{s - (\alpha+n)/p}_{p,p}(\real^d)}

\title[Dyadic norms and Whitney extensions]{Traces of weighted function spaces: dyadic norms and Whitney extensions}

\author[Koskela]{Pekka Koskela}
\address{Department of Mathematics and Statistics, University of Jyv\"askyl\"a, PO~Box~35, FI-40014 Jyv\"askyl\"a, Finland}
\email{pekka.j.koskela@jyu.fi}
\thanks{The authors were supported by the Academy of Finland via the Centre of Excellence in Analysis and Dynamics Research (project No. 307333)}

\author[Soto]{Tom\'as Soto}
\address{Department of Mathematics and Statistics, University of Jyv\"askyl\"a, PO~Box~35, FI-40014 Jyv\"askyl\"a, Finland}
\email{tomas.a.soto@jyu.fi}

\author[Wang]{Zhuang Wang}
\address{Department of Mathematics and Statistics, University of Jyv\"askyl\"a, PO~Box~35, FI-40014 Jyv\"askyl\"a, Finland}
\email{zhuang.z.wang@jyu.fi}



\begin{document}

\maketitle

\begin{abstract}
The trace spaces of Sobolev spaces and related fractional smoothness spaces have been an active area of research since the work of Nikolskii, Aronszajn, Slobodetskii, Babich and Gagliardo among others in the 1950's. In this paper we review the literature concerning such results for a variety of weighted smoothness spaces. For this purpose, we present a characterization of the trace spaces (of fractional order of smoothness), based on integral averages on dyadic cubes, which is well adapted to extending functions using the Whitney extension operator.
\end{abstract}

\section{Introduction}\label{se:introduction}

In 1957, Gagliardo \cite{Ga} gave a characterization of the trace space of the first order Sobolev space $W^{1,p}(\Omega)$, $1<p<\infty$, on a given Lipschitz domain $\Omega \subset \real^d$ in terms of the convergence of a suitable double integral of the boundary values. This work extended the earlier results by Aronszajn \cite{Ar} and Slobodetskii and Babich \cite{SlBa} concerning the case $p=2$. The trace space $B_{p}^{1 - 1/p}(\partial\Omega)$, consisting of all $(d-1)$-Hausdorff measurable functions $u$ on $\partial\Omega$ with
\beqla{eq:diagonal-besov}
  \|u\|_{L^p(\partial\Omega,\mathcal{H}_{d-1})}^p + \int_{\partial \Omega} \int_{\partial \Omega} \frac{|u(x)-u(y)|^p}{|x-y|^{(d-1) + (1-1/p)p}} d\mathcal{H}_{d-1}(x) d\mathcal{H}_{d-1}(y) < \infty,
\eeq
is nowadays commonly called a fractional Sobolev space, a Slobodetskii space or a Besov space. Actually, Gagliardo also verified that the trace space of $W^{1,1}(\Omega)$ is $L^1(\partial \Omega)$ (see also \cite{Mir} for a different proof of this fact). The norm estimates for the trace functions were obtained via Hardy inequalities, while the extension from the boundary was based on a Poisson-type convolution procedure. We refer to the seminal monographs by Peetre \cite{Pe} and Triebel \cite{T} for extensive treatments of the Besov spaces and related smoothness spaces.

A natural variant of this problem asks for the trace spaces associated to weights. Already in 1953, Nikolskii \cite{Ni} had considered the trace problem for Sobolev spaces (for $p = 2$) with weights of the form $x \mapsto \dist(x,\partial \Omega)^{\alpha}$, where $-1 < \alpha < 1$. Other early related results were given by Lizorkin \cite{Li} and Va\v{s}arin \cite{Va}; see \cite{Ne} and \cite{MirRus} for further references. More recently, Tyulenev \cite{Tyu1,Tyu4,Tyu2,Tyu3} has identified the traces of Sobolev and Besov spaces associated to more general Muckenhoupt $A_p$-weights. For related results concerning the traces of weighted Orlicz-Sobolev spaces, we refer to \cite{Fou,Lac,Pal,DhKa,DhKa2} and the references therein.

On the other hand, a notable amount of recent research has focused on extending the theory of Sobolev spaces and related fractional smoothness spaces to the setting of metric measure spaces (including fractal subsets of Euclidean spaces); see e.g. \cite{HKST} and the references therein as well as \cite{HMY}. Works focusing on trace theorems for fractals and related subsets of a Euclidean space include \cite{JoWa, Sa1, Sa2, T4, Shva, IV, Jo, CaHa, HaMa} (we also refer to \cite{TyuVod} for a recent result concerning traces on non-regular subsets of $\real^d$), while trace theorems in more general metric settings have been considered e.g.~in \cite{GKS,SS,LS,MSS,Ma}. In fact, the characterizations of fractional smoothness spaces as retracts of certain sequence spaces in \cite{FJ}, \cite[Section 7]{HMY} and \cite[Proposition 6.3]{BSS} can also be seen as abstract trace theorems.

Motivated by these works, we revisit the Euclidean setting, viewing the upper half-space $\real^{d+1}_+ := \real^d\times(0,\infty)$ as a particularly nice metric space endowed with a weighted measure. We shall introduce equivalent norms for the Besov spaces based integral averages on dyadic cubes. These norms are well adapted for studying the extension of functions defined on $\real^d$ to $\real^{d+1}_+$ via the natural \emph{Whitney extension}. In contrast, the extension operator e.g.~in \cite{MirRus} is based on the Poisson kernel.

Let us begin with a concrete example. We consider functions defined on the real line, but as we will later see, the discussion below generalizes to the setting of higher dimensions as well.

Given $u \in \locint(\real)$ and an interval $I\subset \real$, set
\[
  u(I) := \frac{1}{|I|} \int_I u(x) dx,
\]
where $|I|$ is the length of the interval $I$. For each $k \in \nanu_0$, fix a dyadic decomposition of $\real$ into closed intervals $\{I_{k,i}\}_{i \in \integer}$ so that each $I_{k,i}$ has length $2^{-k}$ and $I_{k,i}\cap I_{k,j}\neq \emptyset$ exactly when $|i-j|\le 1.$  Consider 
the condition
\beqla{eq:eka}
  \|u\|_{L^2(\R)}^2+\sum_{k \in \nanu_0} \sum_{i \in \integer}|u(I_{k,i})-u(I_{k,i+1})|^2<\infty.
\eeq
Now write $Q_{k,i}:=I_{k,i}\times [2^{-k},2^{-k+1}]$ for all admissible $k$ and $i$. Then these squares give us a Whitney decomposition of the upper half-plane $\real^2_+$. Pick a partition of unity in $\bigcup_{k,i} Q_{k,i}$ consisting of functions $\varphi_{k,i}\in C^\infty(\real^2_+)$ such that $|\nabla \varphi_{k,i}|\le 5 \cdot 2^{k}$ and the support of $\varphi_{k,i}$ is contained in a $2^{-k-2}$-neighborhood of $Q_{k,i}$. For $u \in \locint(\real)$, define
\beqla{eq:intro-whitney}
  \ext u := \sum_{k \in \nanu_0,\,i \in \integer} u(I_{k,i})\varphi_{k,i}.
\eeq

Given $f \in W^{1,2}(\real^2_+)$, the trace function $u := \res f\colon \real^2\to \complex$, defined by
\[
  \res f (x) = \lim_{r\to 0} \frac{1}{m_2\big(B\big((x,0),r\big)\cap\real^2_+\big)}\int_{B((x,0),r) \cap \real^2_+} f(y) dm_2(y),
\]
where $m_2$ stands for the $2$-dimensional Lebesgue measure, is well-defined pointwise almost everywhere and satisfies the condition \eqref{eq:eka}. Conversely, if $u \in \locint(\real)$ satisfies \eqref{eq:eka}, we have $\ext u \in W^{1,2}(\real^2_+)$ with the expected norm bound and $\res (\ext u) = u$ pointwise almost everywhere.

We conclude that $u \in \locint(\real)$ belongs to the trace space of $W^{1,2}(\real^2_+)$ if and only if \eqref{eq:eka} holds. Hence the condition \eqref{eq:eka} should characterize the space $B^{1/2}_{2}(\real)$. This is indeed the case; a direct proof is given in Subsection \ref{ss:equivalent-norms} of the Appendix.

Let us next consider the following generalized form of the condition \eqref{eq:eka}: 
\beqla{eq:toka}
  \|u\|_{L^p(\R)}^p+\sum_{k \in \nanu_0} 2^{spk}2^{-k}\sum_{i \in \integer}|u(I_{k,i})-u(I_{k,i+1})|^p<\infty,
\eeq
where $1 < p < \infty$ and $0 < s < 1$. Above we saw that the choice $p=2$ and $s=1/2$ yields the trace space of $W^{1,2}(\real^2_+)$. Similarly, it turns out that the condition \eqref{eq:toka} characterizes the trace space of $W^{1,p}(\real^2_+)$ when $s=1-1/p$.
Where does this value of $s$ come from? The so-called \emph{differential dimension} of the space $B^s_p$ over an $n$-dimensional Euclidean space is $s - n/p$, and the same holds for the space $W^{1,p}$ with $1$ in place of $s$; see e.g.~\cite[Section 3.4.1]{T}. Hence the order of smoothness $s$ of the trace space should satisfy
\[
  s - \frac{1}{p} = 1 - \frac{2}{p},
\]
which rewrites as $s = 1 - 1/p$.

Let us now try to extend a function $u \in \locint(\real)$ to a weighted Sobolev space, by requiring that 
\beqla{eq:intro-weighted}
  \int_{\real^2_+}|\ext u(x)|^p \dist\big(x,\real\times\{0\}\big)^{\alpha}dx + \int_{\real^2_+}|\nabla (\ext u)(x)|^p \dist\big(x,\real\times\{0\}\big)^{\alpha}dx<\infty,
\eeq
where $\alpha > -1$ and $\ext u$ is as defined as above. It turns out that when $\alpha \in (-1,p-1)$, the condition \eqref{eq:intro-weighted} is satisfied when $u$ satisfies \eqref{eq:toka} with $s = 1 - (\alpha+1)/p$. On the other hand, since
\[
  \mu_\alpha\big( B\big( (x,0),r\big) \cap \real^2_+ \big) \approx r^{2+\alpha}
\]
for all $x \in \real$ and $r > 0$, where $\mu_\alpha$ is the measure associated to the weight $x \mapsto \dist(x,\real)^\alpha$, we also see that $\alpha+1 = (2+\alpha) -1 $ appears as a local codimension of $\real$ with respect to the metric measure space $(\real^2_+,\mu_\alpha)$. Hence the drop in the order of the derivative from one to the fractional order $s$ is determined by $p$ and this codimension.

Would the condition \eqref{eq:eka} allow us also to extend functions from $\real$ to a higher-dimensional weighted Euclidean space, e.g.~$(\real^3,\mu_\alpha)$? If so, then the correct condition for the parameter $\alpha$ would be $\alpha > -2$ and the role of $2+\alpha$ above should be taken by $3 + \alpha$. We recover $s = 1/2$ when $(\alpha+2)/p=1/2$, which for $p=2$ gives $\alpha=-1$. This indeed works: \eqref{eq:eka} holds exactly when $u$ is in the trace of  $W^{1,2}(\real^3,\mu_{-1})$, where the measure $\mu_{-1}$ is associated to the weight $x \mapsto \dist(x,\real\times\{0\})^{-1}$ in $\real^3$, and in this case $u$ can be extended as a function in $W^{1,2}(\real^3,\mu_{-1})$ with a suitable modification of the Whitney extension operator \eqref{eq:intro-whitney}.

Can we find yet further function spaces whose traces are characterized by the condition \eqref{eq:eka} or the condition \eqref{eq:toka}? Towards this, let us mention that the space characterized by \eqref{eq:toka} coincides with the \emph{diagonal Triebel-Lizorkin space} $\F^s_{p,p}(\real)$. The scale of Triebel-Lizorkin spaces $\F^s_{p,q}(\real^d)$ on the $d$-dimensional Euclidean space, where $1 \leq p < \infty$, $0 < q \leq \infty$ and $0 < s < 1$, is another widely-studied family of fractional smoothness spaces that arise e.g.~as the complex interpolation spaces between $L^p(\real^d)$ and $W^{1,p}(\real^d)$. The discussion above concerning the traces of weighted Sobolev spaces, with suitable modifications for the parameter ranges, turns out to hold for the traces of these function spaces as well. In particular, when $s \in (0,1)$ and $\alpha \in (-1 , sp-1)$, the condition \eqref{eq:toka} with $s - (\alpha+1)/p$ in place of $s$ characterizes the traces of the functions in $\F^s_{p,q}(\real^2_+,\mu_\alpha)$ for any admissible $q$. A similar trace theorem for the scale of Besov spaces $\B^s_{p,q}(\real^2_+,\mu_\alpha)$ is formulated as Theorem \ref{th:besov-trace} below. The precise definitions of these spaces are given in the next section.

Let us summarize the above discussion. The Whitney extension operator $\ext$ extends a Besov space $B^s_{p}(\real) := \B^s_{p,p}(\real)$ with given smoothness $s \in (0,1)$ linearly and continuously to a number of different (weighted) smoothness spaces defined on $\real^{2}_+$, the trace of all of whose equals $\B^s_{p,p}(\real)$. Moreover, given $n \in \nanu$, a suitable variant of the Whitney extension operator $\ext$ gives us
a similar extension from $\B^s_{p,p}(\real)$ to a variety of (weighted) function spaces defined on $\real^{1+n}$; this is discussed in detail in Subsection \ref{ss:higher-codim}.

To finish discussion, let us state our main results more precisely. Given a pair of function spaces $(X,Y)$, we say that they are a \emph{Whitney trace-extension pair} if $X$ is the trace space of $Y$ in the usual sense and the extension from $X$ to $Y$ is obtained using the natural Whitney extension -- this notion is also defined more precisely in {\bf Definition \ref{de:whitney}} below. The measure $\mu_\alpha$ (where $\alpha > - 1$) below stands for the measure on $\real^{d+1}_+$ defined by
\beqla{eq:weight-measure}
  \mu_\alpha(E)=\int_E w_\alpha \, dm_{d+1},
\eeq
where $w_\alpha \colon \real^{d+1}_+ \to (0,\infty)$ is the weight $(x_1,x_2,\cdots,x_{d+1}) \mapsto \min(1,|x_{d+1}|)^{\alpha}$ and $m_{d+1}$ is the standard Lebesgue measure on $\real^{d+1}_+$. Finally, the definitions of the relevant function spaces are given in {\bf Section \ref{se:definitions}} below.

First off, we have the following trace theorem for the first-order Sobolev spaces.
\thm{th:sobolev-trace}
Let $1\leq p<\infty$ and $-1<\alpha<p-1$. Then $\big(\B^{1-(\alpha+1)/p}_{p, p}(\real^d),W^{1, p}(\real^{d+1}_+, \mu_\alpha)\big)$ is a Whitney trace-extension pair.
\ethm

The analogous trace theorem for the Besov scale reads as follows.
\thm{th:besov-trace}
Let $0 < s < 1$, $1 \leq p \leq \infty$, $0 < q \leq \infty$ and $-1 < \alpha < sp -1$. Then $\big(\bpt,\bpa\big)$ is a Whitney trace-extension pair.
\ethm

Finally, the trace theorem for the Triebel-Lizorkin spaces reads as follows.
\thm{th:triebel-trace}
Let $0 < s < 1$, $1 \leq p < \infty$, $0 < q \leq \infty$ and $-1 < \alpha < sp -1$. Then $\big(\bptd,\fpa\big)$ is a Whitney trace-extension pair.
\ethm

We present a refinement of the case $p = 1$ of Theorem \ref{th:sobolev-trace}, where the Sobolev space $W^{1,1}(\real^{d+1}_+,\mu_\alpha)$ is replaced by a \emph{Hardy-Sobolev space} $h^{1,1}(\real^{d+1}_+,\mu_\alpha)$, in {\bf Section \ref{se:hardy-sobolev}}. The variants of the results above with higher Euclidean codimension are given in {\bf Subsection \ref{ss:higher-codim}} of the Appendix.

The paper is organized as follows. In Section \ref{se:definitions} we give the definitions relevant to our main results and recall some basic properties of the spaces and measures in question. Sections \ref{se:sobolev} through \ref{se:hardy-sobolev} contain the proofs of the aforementioned trace theorems. The Appendix (Section \ref{se:appendix}) deals with various technicalities that we saw fit to postpone from the other sections.

\section{Definitions and preliminaries}\label{se:definitions}

In this section we present the definitions of the relevant function spaces and the Whithey extension operator. Before this, let us introduce some notation that will be used throughout the paper.

\begin{notation}
(i) The majority of this paper will deal with extensions of functions defined on the Euclidean space $\real^d$ to the half-space $\real^{d+1}_+ := \real^d \times (0,\infty)$. The dimension $d \in \nanu := \{1,2,3,\cdots\}$ will be fixed throughout the paper. The $d$-dimensional Lebesgue measure will be denoted by $m_d$. When talking about measures $\mu$ defined on $\real^{d+1}_+$, we may abuse notation by writing $\mu( B(x,r) )$ for $\mu( B(x,r) \cap \real^{d+1}_+)$ when e.g.~$x \in \real^d\times\{0\}$.

\smallskip
(ii) If $(X,\mu)$ is a measure space and $A$ is a $\mu$-measurable subset of $X$ with $0 < \mu(A) < \infty$, we shall write
\[
  f_{A,\mu} := \dashint_A f d\mu := \frac{1}{\mu(A)} \int_A f d\mu
\]
whenever the latter quantity is well-defined, i.e.~when $f \in L^1(A,\mu)$ or $f(x) \geq 0$ for $\mu$-almost every $x \in A$. We may omit $\mu$ from the notation and simply write $f_A$ when $\mu$ is the Lebesgue measure on an Euclidean space and there is no risk of confusion.

\smallskip
(iii) While $\locint(\real^d)$ stands for the space of (complex-valued) locally integrable functions on $\real^d$ in the usual sense, we use the notation $\locint(\real^{d+1}_+)$ with a slightly different meaning: it refers to the space functions that are integrable on bounded subsets of $\real^{d+1}_+$.

\smallskip
(iv) If $f$ and $g$ are two non-negative functions on the same domain, we may use the notation $f \lesssim g$ with the meaning that $f \leq C g$ in the domain, where the constant $C > 0$ is usually independent of some parameters obvious from the context. The notation $f \approx g$ means that $f \lesssim g$ and $g \lesssim f$.
\end{notation}

\defin{de:sobolev}
Suppose that $\mu$ is a Borel regular measure on $\real^d$ such that every Euclidean ball has positive and finite $\mu$-measure. 

Let $p \in [1,\infty)$. Then $W^{1, p}(\real^d, \mu)$ is defined as the normed space of measurable functions $f \in \locint(\real^d)$ such that the first-order distributional derivatives of $f$ coincide with functions in $\locint(\real^d)$ and
\begin{equation}\label{eq:sobolev}
  \|f\|_{W^{1, p}(\real^d, \mu)}:= \|f\|_{L^p(\real^d,\mu)} +  \|\nabla f\|_{L^p(\real^d,\mu)}
\end{equation}
is finite.

The space $W^{1, p}(\real^{d+1}, \mu)$ is defined similarly, by replacing $\real^d$ with $\real^{d+1}_+$ in \eqref{eq:sobolev}.
\edefin

In order to formulate the dyadic norms of the relevant fractional smoothness spaces, we recall the standard dyadic decompositions of $\real^d$ and $\real^{d+1}_+$. Denote by $\dyadic_d$ the collection of dyadic semi-open cubes in $\real^d$, i.e.~the cubes of the form $Q := 2^{-k}\big((0,1]^d + m\big)$, where $k \in \integer$ and $m \in \integer^d$, and $\dyadic^+_d$ for the cubes in $\dyadic_d$ which are contained in the upper half-space $\real^{d-1}\times(0,\infty)$. Write $\ell(Q)$ for the edge length of $Q \in \dyadic_d$, i.e.~$2^{-k}$ in the preceding representation, and $\dyadic_{d,k}$ for the cubes $Q \in \dyadic_d$ such that $\ell(Q) = 2^{-k}$. If $x \in \real^d$ (resp.~$x \in \real^{d+1}_+$) and $k \in \integer$, we may write write $Q^x_k$ for the unique cube in $\dyadic_d$ (resp.~$\dyadic^+_{d+1}$) such that $x \in Q$ and $\ell(Q) = 2^{-k}$.

We say that $Q$ and $Q'$ in $\dyadic_d$ are neighbors and write $Q \sim Q'$ if $\frac12 \leq \ell(Q)/\ell(Q') \leq 2$ and $\overline{Q}\cap\overline{Q'} \neq \emptyset$. Note that every $Q$ has a uniformly finite number of neighbors.

\defin{de:besov} Suppose that $\mu$ is a Borel regular measure on $\real^d$ such that every Euclidean ball has positive and finite $\mu$-measure.

Let $s \in (0,1)$, $p \in [1,\infty]$ and $q \in (0,\infty]$. Then the Besov space $\bp$ is defined as the normed (or quasi-normed when $q < 1$) space of functions $f \in \locint(\real^d)$ such that 
\begin{equation}\label{eq:besov-norm}
  \|f\|_{\bp} := \|f\|_{L^p(\real^d,\mu)} + \bigg( \sum_{k = 0}^\infty 2^{ksq}\Big( \sum_{Q \in \dyadic_{d,k}}\mu(Q) \sum_{Q' \sim Q} \big|f_{Q,\mu} - f_{Q',\mu}\big|^p \Big)^{q/p} \bigg)^{1/q}
\end{equation}
(standard modification for $p = \infty$ and/or $q = \infty$) is finite.
\edefin

\defin{de:triebel} Suppose that $\mu$ is a Borel regular measure on $\real^d$ such that every Euclidean ball has positive and finite $\mu$-measure.

Let $s \in (0,1)$, $p \in [1,\infty)$ and $q \in (0,\infty]$. Then the Triebel-Lizorkin space $\fp$ is defined as the normed (or quasi-normed when $q < 1$) space of functions $f \in \locint(\real^d)$ such that 
\begin{equation}\label{eq:triebel-norm}
  \|f\|_{\fp} := \|f\|_{L^p(\real^d,\mu)} + \bigg(\int_{\real^d} \Big(\sum_{k = 0}^\infty 2^{ksq} \sum_{Q' \sim Q^x_k} \big|f_{Q^x_k,\mu} - f_{Q',\mu}\big|^q \Big)^{p/q}\, d\mu(x) \bigg)^{1/p}
\end{equation}
(standard modification for $q = \infty$) is finite.
\edefin

The spaces $\B^s_{p,q}(\real^{d+1}_+,\mu)$ and $\F^s_{p,q}(\real^{d+1}_+,\mu)$ are defined similarly, by replacing $\real^d$ with $\real^{d+1}_+ := \real^d\times(0,\infty)$ in \eqref{eq:besov-norm} and \eqref{eq:triebel-norm} respectively, and omitting the terms corresponding to the cubes $Q \in \dyadic_{d+1}\setminus \dyadic^+_{d+1}$ and $Q' \in \dyadic_{d+1}\setminus\dyadic^+_{d+1}$.

\begin{remark}\label{re:basic-properties}
One routinely checks that $\bp$ and $\fp$ are quasi-Banach spaces (Banach spaces for $q \geq 1$). Fubini's theorem implies that
\[
  \F^{s}_{p,p}(\real^d,\mu) = \B^{s}_{p,p}(\real^d,\mu)
\]
with equivalent norms for $p \in [1,\infty)$, and the monotonicity of the $\ell^q$-norms shows that
\[
  \B^{s}_{p,q}(\real^d,\mu) \subset \B^{s}_{p,q'}(\real^d,\mu) \quad \text{and} \quad \F^{s}_{p,q}(\real^d,\mu) \subset \F^{s}_{p,q'}(\real^d,\mu)
\]
with continuous embeddings when $q' > q$. All this of course holds with $\real^{d+1}_+$ in place of $\real^d$.
\end{remark}

In case $\mu$ is the standard Lebesgue measure on $\real^d$, we shall omit $\mu$ from the notation of the three function spaces above and simply write $W^{1,p}(\real^d)$, $\bpd$ and $\fpd$ where appropriate.

\begin{remark}
(i) A Besov quasinorm that is perhaps more standard in the literature is given by
\begin{equation}\label{standard-besov}
  f \mapsto \|f\|_{L^p(\real^d,\mu)} + \bigg( \int_{0}^\infty t^{-sq }\Big( \int_{\real^d} \dashint_{B(x,t)} |f(x) - f(y)|^p d\mu(y) d\mu(x) \Big)^{q/p} \frac{dt}{t} \bigg)^{1/q}.
\end{equation}
A straightforward calculation using Fubini's theorem shows that if $q = p$ and $\mu = m_d$, then then the $p$th power of this this quasinorm is comparable to
\[
  \|f\|_{L^p(\real^d)}^p + \int_{\real^d}\int_{\real^d} \frac{|f(x)-f(y)|^p}{|x-y|^{d+sp}} dx dy,
\]
which is of the same form as the quantity \eqref{eq:diagonal-besov} in the introduction.

\smallskip
(ii) To justify the Definitions \ref{de:besov} and \ref{de:triebel} above, let us point out that if the measure $\mu$ is doubling with respect to the Euclidean metric, i.e.~if there exists a constant $c \geq 1$ such that
\[
  \mu\big( B(x,2r) \big) \leq c \mu\big( B(x,r)\big) \quad \text{ for all } x \in \real^d\text{ and } r > 0,
\]
then the quasi-norm \eqref{eq:besov-norm} is comparable to the quasi-norm \eqref{standard-besov} above. We refer to Subsection \ref{ss:equivalent-norms} of the Appendix for details.

\smallskip
(iii) Quasinorms similar to \eqref{eq:besov-norm} and \eqref{eq:triebel-norm} in the setting of metric measure spaces were also considered in \cite[Definition 5.1]{SS} in terms of a \emph{hyperbolic filling} of $\real^d$. Another similar variant in the weighted Euclidean setting has been considered in \cite{Tyu3}.
\end{remark}

We now give the definitions corresponding to the Whitney extensions discussed in the introduction. To this end, we have to define a partition of unity corresponding to the standard Whitney decomposition of the half-space $\real^{d+1}_+$. For $Q \in \dyadic_{d,k}$, $k \in \integer$, write $\W(Q) := Q \times (2^{-k},2^{-k+1}] \in \dyadic^+_{d+1,k}$. To simplify the notation in the sequel, we further define $\dyadic^0_d := \cup_{k \geq 0} \dyadic_{d,k}$.

It is then easy to see that $\{\W(Q) \,:\, Q \in \dyadic_d\}$ is a Whitney decomposition of $\real^{d}\times(0,\infty)$ with respect to the boundary $\real^d\times\{0\}$. For all $Q \in \dyadic^0_d$, define a smooth function $\psi_Q \colon \real^{d+1}_{+}\to[0,1]$ such that $\Lip \psi_Q \lesssim 1/\ell(Q)$, $\inf_{x \in \W(Q)} \psi_Q(x) > 0$ uniformly in $Q \in \dyadic^0_d$, $\supp \psi_Q$ is contained in an $\frac{\ell(Q)}{4}$-neighborhood of $\W(Q)$ and
\[
  \sum_{Q \in \dyadic^0_d} \psi_Q \equiv 1 \quad \text{in} \quad \bigcup_{Q \in \dyadic^0_d} \W(Q).
\]
Let us point out that the sum above is locally finite -- more precisely, it follows from the definition that
\beqla{eq:bump-support}
  \supp \psi_Q \cap \supp \psi_{Q'} \neq \emptyset \quad \text{if and only if} \quad Q \sim Q'.
\eeq

\defin{de:whitney}
(i) Let $f \in \locint(\real^d)$. Then the Whitney extension $\ext f \colon \real^{d+1}_+\to\complex$ is defined by
\[
  \ext f(x) = \sum_{Q \in \dyadic^0_d} \Big( \dashint_{Q} f dm_d \Big) \psi_Q(x).
\]
This definition gives rise in the obvious way to the linear operator $\ext \colon \locint(\real^d) \to C^{\infty}(\real^{d+1}_+)$.

(ii) Let $X \subset \locint(\real^d)$ be a quasinormed function space on $\real^d$, and let $Y$ be a quasinormed function space on the weighted half-space $(\real^{d+1}_+,\mu)$. We say that $(X,Y)$ is a \emph{Whitney trace-extension pair} if $\ext$ maps $X$ continuously into $Y$, if the trace function $\res f$ defined by
\beqla{eq:lebesgue}
  \res f(x) = \lim_{r\to 0} \dashint_{B((x,0),r)\cap \real^{d+1}_+} f(y) d\mu(y),
\eeq
is for all $f \in Y$ well-defined almost everywhere and belongs to $\locint(\real^d)$, if $\res$ maps $Y$ continuously into $X$ and if
\[
  \res ( \ext f ) = f
\]
pointwise almost everywhere for all $f \in X$.
\edefin

For the proofs of our main results, let us recall some basic facts about the weights $w_\alpha$ and measures $\mu_\alpha$ defined in \eqref{eq:weight-measure}. First, it is well-known that for $\alpha > -1$, the weight $w_\alpha$ belongs to the Muckenhoupt class $A_r$ for all $r > \max(\alpha+1,1)$, which implies that the measure $\mu_\alpha$ satisfies the doubling property with respect to the standard Euclidean metric (see e.g.~\cite[Chapter 15]{HKM} or \cite{DILTV}). This in particular means that
\[
  \mu_\alpha\big( Q \big) \approx \mu_\alpha\big( Q' \big) \quad \text{if} \quad Q \sim Q'.
\]
A straightforward calculation also shows that
\beqla{eq:balls-at-boundary}
  \mu_\alpha\big( B(x,r) \big) \approx r^{d+1+\alpha}
\eeq
for all $x \in \real^d\times\{0\}$ and $0 < r \leq 1$.

Finally, let us recall the standard $(1,1)$-Poincar\'e inequality satisfied by the functions that are locally $W^{1,1}$-regular in the upper half-space. If $Q$ is a cube in $\real^{d+1}_+$ such that $\dist(Q,\real^d\times\{0\}) > 0$ and $f \in W^{1,1}(Q)$, we have
\beqla{poincare}
  \dashint_Q |f-f_{Q}|\, dm_{d+1}\leq C \ell(Q)\dashint_Q |\nabla f| dm_{d+1}
\eeq
for some constant $C$ independent of $Q$ and $f$.

\section{Proof of Theorem \ref{th:sobolev-trace}}\label{se:sobolev}
\begin{proof}
(i) Let us first prove the desired norm inequality for the Whitney extension $\ext f$ of $f \in \B^{1-(\alpha+1)/p}_{p,p}(\real^d)$. We begin by noting that if $Q \in \dyadic^0_d$, it follows directly from the definitions that $w_\alpha \approx \ell(Q)^\alpha$ in $\W(Q)$, and hence
we have $\mu_\alpha(\W(Q)) \approx \ell(Q)^{\alpha} m_{d+1}(\W(Q)) \approx \ell(Q)^{d+1+\alpha}$. Since the supports of the functions $\psi_Q$ have bounded overlap, the $L^p(\real^{d+1}_+,\mu_\alpha)$-norm of $\ext f$ is thus easy to estimate:
\begin{align}
  \int_{\real^{d+1}_+} |\ext f|^p d\mu_\alpha & \lesssim \sum_{Q \in \dyadic^0_d} \mu_\alpha\big(\W(Q)\big) \dashint_Q |f|^p dm_d \approx \sum_{Q \in \dyadic^0_d} \ell(Q)^{\alpha+1} \int_Q |f|^p dm_d \notag \\
& = \sum_{k \geq 0} 2^{-k(\alpha+1)} \sum_{Q \in \dyadic_{d,k}} \int_Q |f|^p dm_d =\sum_{k \geq 0} 2^{-k(\alpha+1)} \int_{\real^d}|f|^p\, dm_d \notag \\
&=\sum_{k \geq 0} 2^{-k(\alpha+1)}  \|f\|_{L^p(\real^d)}^p\approx \|f\|_{L^p(\real^d)}^p.\label{eq:sobolev-lp-estimate}
\end{align}

In order to estimate the $L^p(\real^{d+1}_+,\mu_\alpha)$-norm of $|\nabla(\ext f)|$, we divide the half-space $\real^{d+1}_+$ into two parts: $X_1:=\bigcup_{P\in \dyadic_d^0}\W(P)$ and $X_2:=\real^{d+1}_+\setminus X_1$. Now if $x \in X_1$, i.e.~$x\in \W(P)$ for some $P \in \dyadic^0_d$, we have that $\sum_{Q\in\dyadic_d^0}\psi_Q(x)=1$, and as noted in \eqref{eq:bump-support}, the terms in this sum are nonzero at most for the cubes $Q$ such that $Q \sim P$. Hence
\begin{align*}
  \ext f(x)-\dashint_{P} f\, dm_d&=\sum_{Q\in \dyadic_d^0}\Big(\dashint_Q f\, dm_d\Big)\psi_Q(x)-\dashint_{P} f\, dm_d\\
&=\sum_{Q\sim P}\bigg(\dashint_Q f\, dm_d-\dashint_{P} f\, dm_d\bigg)\psi_Q(x)=\sum_{Q\sim P}\big(f_Q-f_{P}\big)\psi_Q(x),
\end{align*}
and the Lipschitz continuity of the functions $\psi_Q$ yields
\begin{align}
|\nabla(\mathcal Ef)(x)|&\leq |\Lip(\mathcal Ef)(x)|=\Big|\Lip\Big(\ext f(\cdot)-\dashint_{P} f\, dm_d\Big)(x)\Big| \notag \\
&\leq \sum_{Q\sim P}|f_Q-f_{P}||\Lip(\psi_Q)(x)|\lesssim \sum_{Q\sim P}\frac{1}{\ell(Q)}|f_Q-f_{P}|. \label{eq:common-energy-estimate-1}
\end{align}
This means that
\begin{align}
  \int_{X_1} |\nabla(\ext f)|^p d\mu_{\alpha} & = \sum_{P \in \dyadic^0_d} \int_{\W(P)} |\nabla(\ext f)|^p d\mu_{\alpha} \lesssim \sum_{P \in \dyadic^0_d} \mu_\alpha\big(\W(P)\big) \sum_{Q\sim P}\frac{1}{\ell(Q)^p}|f_Q-f_{P}|^p \notag\\
& \approx \sum_{P \in \dyadic^0_d} \ell(P)^{d+1+\alpha} \sum_{Q\sim P}\frac{1}{\ell(Q)^p}|f_Q-f_{P}|^p \notag \\
& \approx \sum_{P \in \dyadic^0_d} \ell(P)^{-(1 - \frac{\alpha+1}p)p} m_d(P) \sum_{Q\sim P}|f_Q-f_{P}|^p \notag\\
& \lesssim \|f\|_{\B^{1-(\alpha+1)/p}_{p,p}(\real^d)}^p. \label{eq:sobolev-energy-estimate-1}
\end{align}

If on the other hand $x\in X_2$, we can have $\psi_Q(x) \neq 0$ only for $Q \in \dyadic_{d,0}$. Thus,
\[
  \ext f(x) = \sum_{Q \in \dyadic_{d, 0} }\Big(\dashint_{Q} f dm_d\Big) \psi_Q(x)=\sum_{\substack{Q \in \dyadic_{d,0}\\\supp \psi_Q \owns x }}f_Q\psi_Q(x),
\]
and using the Lipschitz continuity of the functions $\psi_Q$ as above, we get
\[
  |\nabla(\mathcal Ef)(x)|\leq |\Lip(\mathcal Ef)(x)| \leq \sum_{\substack{Q \in \dyadic_{d,0}\\\supp \psi_Q \owns x}} |f_Q||\Lip(\psi_Q)(x)|\lesssim \sum_{Q \in \dyadic_{d,0}}|f_Q| \chi_{\supp \psi_Q}(x).
\]
Since $\mu_{\alpha} (\supp \psi_Q) \approx \mu_{\alpha} (\W(Q)) \approx 1$ for all $Q \in \dyadic_{d,0}$, the estimate above yields
\beqla{eq:sobolev-energy-estimate-2}
  \int_{X_2} |\nabla(\ext f)|^p d\mu_{\alpha} \lesssim \sum_{Q \in \dyadic_{d,0}} |f_Q|^p \leq \sum_{Q \in \dyadic_{d,0}} \int_Q |f|^p dm_d = \|f\|_{L^p(\real^d)}^p.
\eeq
Combining \eqref{eq:sobolev-lp-estimate}, \eqref{eq:sobolev-energy-estimate-1} and \eqref{eq:sobolev-energy-estimate-2}, we arrive at
\[
  \|\ext f\|_{L^p(\real^{d+1_+},\mu_\alpha)} + \| \nabla(\ext f) \|_{L^p(\real^{d+1_+},\mu_\alpha)} \lesssim \|f\|_{\B^{1-(\alpha+1)/p}_{p,p}(\real^d)},
\]
which is the desired norm inequality.

\smallskip
(ii) Let us now consider the existence and norm of the trace function $\res f$ of a function $f\in W^{1, p}(\real^{d+1}_+, \mu_\alpha)$. For $k \in \nanu_0$, define the function $\T_k f \colon \real^d\to\complex$ by
\[
  \T_k f := \sum_{Q \in \dyadic_{d,k}} \Big(\dashint_{\N(Q)} f dm_{d+1} \Big) \chi_Q,
\]
where $\N(Q) := \frac54\W(Q) := \{y\in \real^{d+1}_+ \,:\, \dist(y,\W(Q)) < \frac14 \ell(Q) \}$ -- note that the functions $\T_k f$ are well-defined, since $f \in L^1(\N(Q),\mu_\alpha)$ implies $f \in L^1(\N(Q),m_{d+1})$ for all $Q \in \dyadic^0_d$. We first show that the limit $\lim_{k\to\infty} \T_k f$ exists pointwise $m_d$-almost everywhere in $\real^d$ (and, in fact, in $L^p(\real^d)$). The limit function will be called $\res f$ for now even though it is not of the same form as in Definition \ref{de:whitney} -- we shall return to this point in part (iii) below.

To verify the existence of the limit in question, it suffices to show that the function
\[
  f^* := \sum_{k \geq 0} \big|\T_{k+1}f - \T_k f \big| + \big|\T_0 f\big|
\]
belongs to $L^p(\real^d)$.

Let $P \in \dyadic_{d,0}$. Because $m_{d+1}(\N(P)) \approx 1$ and $\mu_\alpha \approx 1$ in $\N(P)$, we get
\begin{align}
  \int_P |f^*|^p dm_d & \lesssim \int_P \Big| \sum_{k \geq 0} \Big( \T_{k+1} f(x) - \T_k f(x) \Big) \Big|^p dm_d (x) + \int_{\N(P)} |f|^p dm_{d+1} \notag \\
  & \approx \int_P \Big| \sum_{k \geq 0} \Big( \T_{k+1} f(x) - \T_k f(x) \Big) \Big|^p dm_d (x) + \int_{\N(P)} |f|^p d \mu_\alpha \notag \\
  & = \int_P  \Big(\sum_{k \geq 0} 2^{-k\epsilon/p}\big|2^{k\epsilon/p}\big( \T_{k+1} f(x) - \T_k f(x) \big) \big|\Big)^p dm_d (x) + \int_{\N(P)} |f|^p d \mu_\alpha \notag \\
  & \lesssim \sum_{k \geq 0} 2^{k\epsilon} \int_{P} \big| \T_{k+1} f(x) - \T_k f(x) \big|^p dm_d(x) + \int_{\N(P)} |f|^p d \mu_\alpha, \label{eq:lp-estimate-s}
\end{align}
where $\epsilon := p - (\alpha + 1) > 0$ and the last estimate uses H\"older's inequality.

In order to estimate the $k$th integral above, recall that for $x \in \real^d$, $Q_k^x$ stands for unique cube in $\dyadic_{d, k}$ that contains $x$. By the definition of the $\N(Q)$'s, the intersection of $\N(Q_k^x)$ and $\N(Q_{k+1}^x)$ contains a cube $\tilde Q$ with edge length comparable to $2^{-k}$. We thus have the estimate
\begin{align*}
|\T_kf(x)-\T_{k+1}f(x)|&=\bigg|\dashint_{\N(Q_k^x)} f\, dm_{d+1}-\dashint_{\N(Q_{k+1}^x)} f\, dm_{d+1}\bigg|\\
&\leq \bigg|\dashint_{\N(Q_k^x)} f\, dm_{d+1}-\dashint_{\tilde Q} f\, dm_{d+1}\bigg|\\
& \qquad \qquad \qquad +\bigg|\dashint_{\tilde Q} f\, dm_{d+1}-\dashint_{\N(Q_{k+1}^x)} f\, dm_{d+1}\bigg|\\
&\lesssim \dashint_{\N(Q_k^x)}|f-f_{\N(Q_k^x)}|\, dm_{d+1} +\dashint_{\N(Q_{k+1}^x)}|f-f_{\N(Q_{k+1}^x) }|\, dm_{d+1}.
\end{align*}
We have $w_\alpha(y) \approx 2^{-k\alpha}$ for all $y \in \N(Q_k^x)$, and hence also $\mu_\alpha(\N(Q_k^x)) \approx 2^{-k\alpha} m_{d+1}(\N(Q_k^x))$ as in part (i) above. We may therefore use the Poincar\'e inequality \eqref{poincare} in conjunction with H\"older's inequality to estimate the first integral from above by
\begin{align*}
2^{-k}\dashint_{\N(Q_k^x)}|\nabla f|\, dm_{d+1}\approx 2^{-k}\dashint_{\N(Q_k^x)}|\nabla f|\, d\mu_\alpha \leq 2^{-k}\bigg(\dashint_{\N(Q_k^x)} |\nabla f|^pd\mu_\alpha\bigg)^{1/p}.
\end{align*}
A similar estimate obviously holds for the second integral. We thus get
\begin{equation}\label{poincare-estimate}
|\T_kf(x)-\T_{k+1}f(x)|\lesssim2^{-k}\bigg(\dashint_{\N(Q_k^x)} |\nabla f|^pd\mu_\alpha\bigg)^{1/p}+2^{-k}\bigg(\dashint_{\N(Q_{k+1}^x)} |\nabla f|^pd\mu_\alpha\bigg)^{1/p},
\end{equation}
and hence
\begin{align}
  &\int_{P} \big| \T_{k+1} f(x) - \T_k f(x) \big|^p dm_d(x) 
   = \sum_{\substack{Q \in \dyadic_{d,k}\\Q\subset P}} \int_{Q} \big| \T_{k+1} f(x) - \T_k f(x) \big|^p dm_d(x) \notag\\
 \lesssim &\sum_{\substack{Q \in \dyadic_{d,k}\\Q\subset P}} m_d(Q) \sum_{\substack{Q' \in \dyadic_{d,k+1}\\ Q'\subset Q}}  \bigg(\ell(Q)^p\dashint_{\N(Q)}|\nabla f|^p\, d\mu_\alpha+\ell(Q')^p\dashint_{\N(Q')}|\nabla f|^p\, d\mu_\alpha\bigg) \notag\\
 \lesssim &\ \  2^{-k(d+p)}\sum_{\substack{Q \in \dyadic_{d,k}\cup \dyadic_{d,k+1}\\Q\subset P}} \dashint_{\N(Q)}|\nabla f|^p\, d\mu_\alpha=2^{-k\epsilon}\sum_{\substack{Q \in \dyadic_{d,k}\cup \dyadic_{d,k+1}\\Q\subset P}} \int_{\N(Q)}|\nabla f|^p\, d\mu_\alpha.\notag
\end{align}
Plugging this into \eqref{eq:lp-estimate-s} and summing over $P \in \dyadic_{d,0}$, we arrive at
\begin{align*}
\|f^*\|_{L^p(\real^d,m_d)}^p&\lesssim  \sum_{P\in \dyadic_{d, 0}} \sum_{k \geq 0} \sum_{\substack{Q \in \dyadic_{d,k}\cup \dyadic_{d,k+1}\\Q\subset P}} \int_{\N(Q)}|\nabla f|^p\, d\mu_\alpha +\sum_{P\in \dyadic_{d, 0}}\int_{\N(P)} |f|^p d \mu_\alpha\\
& \approx \sum_{Q\in \dyadic_{d}^0}\int_{\N(Q)}|\nabla f|^p\, d\mu_\alpha +\sum_{P\in \dyadic_{d, 0}}\int_{\N(P)} |f|^p d \mu_\alpha\\
& \lesssim \|f\|^p_{W^{1, p}(\real^{d+1}, \mu_\alpha)}.
\end{align*}
Here the last inequality follows from the fact that $\sum_{Q\in\dyadic_d^0}\chi_{\N(Q)}\leq 2$.

Hence $f^*(x) < \infty$ for $m_d$-almost every $x \in \real^d$, so the limit $\res f(x) := \lim_{k\to\infty} \T_k f(x)$ exists at these points. In the remainder of this proof, we shall abuse notation by writing simply $f$ for $\res f$. Since $|f| \leq |f^*|$ almost everywhere in $\real^d$, the estimate above immediately gives
\[
  \|f\|_{L^p(\real^d)} \lesssim \|f\|_{W^{1, p}(\real^{d+1}, \mu_\alpha)}.
\]

Now to estimate the $\B^{1-(1+\alpha)/p}_{p, p}$-energy of $f$, let $Q \in \dyadic_{d,k}$ with $k \geq 0$ and write $Q^* := Q \cup \bigcup_{Q' \sim Q} Q'$. We get
\begin{align*}
  \sum_{Q' \sim Q} \big|f_{Q} - f_{Q'}\big|^p & \lesssim \sum_{Q' \sim Q} \Big( \big|f_{Q} - f_{\N(Q)}\big|^p + \big|f_{Q'} - f_{\N(Q')}\big|^p + \big|f_{\N(Q)} - f_{\N(Q')}\big|^p\Big) \\
& \lesssim \dashint_{Q^*} \big| f(x) - \T_k f(x) \big|^p dm_d(x) + \sum_{Q' \sim Q} \big|f_{\N(Q)} - f_{\N(Q')}\big|^p.
\end{align*}
Note that $m_d(Q^*) \approx m_d(Q)$, that the collection of cubes $\{ Q^* \,:\, Q \in \dyadic_{d,k} \}$ has bounded overlap (uniformly in $k$) and that $m_d(Q) / \mu_\alpha(\N(Q)) \approx 2^{k(\alpha+1)}$. Using these facts together with an estimate similar to \eqref{poincare-estimate}, we get
\begin{align}
  & \sum_{Q \in \dyadic_{d,k}} m_d(Q) \sum_{Q' \sim Q} \big|f_{Q} - f_{Q'}\big|^p \notag\\
\lesssim & \int_{\real^d} \big| f(x) - \T_k f(x) \big|^p dm_d(x) \notag\\
& \qquad \qquad + 2^{k(\alpha+1)} \sum_{Q \in \dyadic_{d,k}} \mu_\alpha\big(\N(Q)\big)\sum_{Q' \sim Q} \big|f_{\N(Q)} - f_{\N(Q')}\big|^p \label{eq:common-energy-estimate-2} \\
\lesssim & \int_{\real^d} \big| f(x) - \T_k f(x) \big|^p dm_d(x) + 2^{k(\alpha+1 - p)} \int_{\cup_{2^{-k-1} \leq \ell(Q') \leq 2^{-k+1}} \N(Q')} |\nabla f|^pd\mu_\alpha \notag \\
& =: I_k + 2^{k(\alpha+1-p)} I'_k, \notag
\end{align}
so that
\beqla{eq:sobolev-trace-estimate-1}
  \sum_{k \geq 0} 2^{k(1-\frac{\alpha+1}{p})p} \sum_{Q \in \dyadic_{d,k}} m_d(Q) \sum_{Q' \sim Q} \big|f_{Q} - f_{Q'}\big|^p \lesssim \sum_{k \geq 0} 2^{k(1-\frac{\alpha+1}{p})p} I_k + \sum_{k \geq 0} I'_k.
\eeq
We have
\beqla{eq:sobolev-trace-estimate-2}
  \sum_{k \geq 0} I'_k \lesssim \|f\|_{W^{1,p}(\real^{d+1}_+,\mu_\alpha)}^p
\eeq
because the domains of integration in the definition of the $I'_k$'s have bounded overlap. To estimate the terms $I_k$, we may take $\epsilon \in (0,p-\alpha-1)$ and proceed as in the estimates following \eqref{eq:lp-estimate-s}:
\begin{align*}
  I_k & \lesssim \sum_{n \geq k} 2^{(n-k)\epsilon} \int_{\real^d} \big| \T_{n+1}f(x) - \T_n f(x) \big|^p dm_d(x) \\
& \lesssim \sum_{n \geq k} 2^{(n-k)\epsilon} 2^{-n(d+p)} \sum_{Q \in \dyadic_{d,n}\cup \dyadic_{d,n+1}} \dashint_{\N(Q)}|\nabla f|^p \, d\mu_\alpha\\
& \approx \sum_{n \geq k} 2^{(n-k)\epsilon} 2^{-n(p - \alpha - 1)} \sum_{Q \in \dyadic_{d,n}\cup \dyadic_{d,n+1}} \int_{\N(Q)}|\nabla f|^p  \, d\mu_\alpha \\
& =: \sum_{n \geq k} 2^{(n-k)\epsilon} 2^{-n(p - \alpha - 1)} O'_n,
\end{align*}
so that
\begin{align*}
  \sum_{k \geq 0} 2^{k(1-\frac{\alpha+1}{p})p} I_k & \lesssim \sum_{n \geq 0} 2^{n(\alpha+1-p + \epsilon)} O'_n \sum_{0 \leq k \leq n} 2^{k(p - \alpha - 1 - \epsilon) } \approx \sum_{n \geq 0} O'_n \notag
\\ & \lesssim \|f\|_{W^{1,p}(\real^{d+1}_+,\mu_\alpha)}
\end{align*}
where the last estimate follows from the definition of the norm.
Plugging this and \eqref{eq:sobolev-trace-estimate-2} into \eqref{eq:sobolev-trace-estimate-1}, we get the desired energy estimate for $\res f$.

\smallskip
(iii) Let $\res$ be as in part (ii) above. Since $m_d$-almost all points of $\real^d$ are Lebesgue points of a function $f \in \B^{1-(1+\alpha)/p}_{p, p}$, it is evident from the definition of $\res$ that $\res( \ext f) = f$ pointwise $m_d$-almost everywhere.

We are now done with the proof of the Theorem, with the exception that the trace operator $\res$ considered in part (ii) is not of the form required by Definition \ref{de:whitney}. This is in fact a cosmetic difference -- by a well-known argument, if $f \in W^{1,p}(\real^{d+1}_+,\mu_\alpha)$, then the point $(x,0)$ is for $m_d$-almost all $x \in \real^d$ in a sense a $\mu_\alpha$-Lebesgue point of $f$. We refer to Subsection \ref{ss:lebesgue} for details. Keeping this fact in mind, it is easily seen that the function $\res f$ considered in part (ii) coincides almost everywhere with the function in \eqref{eq:lebesgue} (with $\mu = \mu_\alpha$).
\end{proof}

\section{Proof of Theorem \ref{th:besov-trace}}\label{se:besov}
\begin{proof}
For simplicity, we only consider the case $q = p < \infty$. The cases where $q \in (0,\infty]$ and/or $p = \infty$ can be proven by simple modifications of the arguments below.

\smallskip
(i) We first establish the desired norm inequality for the function $\ext f$ for $f \in \bptd$. To begin with, since the parameters $p$, $s$ and $\alpha$ are also admissible for Theorem \ref{th:sobolev-trace}, the estimate \eqref{eq:sobolev-lp-estimate} therein tells us that
\beqla{eq:lp-estimate-2}
  \| \ext f \|_{L^p(\real^{d+1}_+,\mu_\alpha)} \lesssim \|f\|_{L^p(\real^d)}.
\eeq
Now to estimate the $\bpad$-energy of $\ext f$, we divide the dyadic cubes in $\real^{d+1}_+$ into three classes that will be considered separately. For $k \geq 0$, write $\dyadic^{1}_k$ for the collection of dyadic cubes $Q$ in $\dyadic^+_{d+1}$ with edge length $2^{-k}$ such that $\dist(Q,\real^d\times\{0\}) \geq 2$, $\dyadic^{2}_k$ for the collection of dyadic cubes $Q$ in $\dyadic^+_{d+1}$ with edge length $2^{-k}$ such that $2^{-k} \leq \dist(Q,\real^d\times\{0\}) < 2$ and $\dyadic^{3}_k$ for the collection of dyadic cubes in $\dyadic^+_{d+1}$ with edge length $2^{-k}$ whose closures intersect $\real^d\times\{0\}$. Also write $\dyadic^{2,*}_k$ for the collection of cubes in  $\cup_{i=\max(k-1,0)}^{k+1} \dyadic^{2}_i$ that are contained in $\cup_{Q \in \dyadic^{2}_k} Q$.

We thus want to estimate
\begin{align}
& \sum_{Q \in \dyadic^{1}_k}\mu_\alpha(Q) \sum_{Q' \sim Q} \big|(\ext f)_{Q,\mu_\alpha} - (\ext f)_{Q',\mu_\alpha}\big|^p + \sum_{Q \in \dyadic^{2}_k}\mu_\alpha(Q) \sum_{\substack{Q' \sim Q \\ Q' \in \dyadic^{2,*}_k}} \big|(\ext f)_{Q,\mu_\alpha} - (\ext f)_{Q',\mu_\alpha}\big|^p \notag\\
& \qquad \qquad \qquad + \sum_{Q \in \dyadic^{3}_k}\mu_\alpha(Q) \sum_{Q' \sim Q} \big|(\ext f)_{Q,\mu_\alpha} - (\ext f)_{Q',\mu_\alpha}\big|^p =: O^{(1)}_k + O^{(2)}_k + O^{(3)}_k \label{eq:cubes-distances}
\end{align}
at each level $k \geq 0$ -- the reason why we can omit the terms corresponding to $Q' \notin \dyadic^{2,*}_k$ in the middle sum is that a comparable term is contained in $O^{(1)}_k$, $O^{(3)}_k$, $O^{(2)}_{k+1}$ or $O^{(1)}_{k-1}$.

We first note that $O^{(1)}_k$ can for $k \in \{0,1\}$ be simply estimated by
\[
  \|\ext f\|_{L^p(\real^{d+1}_+,\mu_\alpha)}^p \lesssim \|f\|_{L^p(\real^d)}^p.
\]
Now suppose that $Q \in \dyadic^{1}_k$ with $k \geq 2$ and $Q' \sim Q$. Using the Lipschitz continuity of the bump functions $\psi_P$ and noting that we can only have $\supp \psi_P \cap (Q\cup Q') \neq \emptyset$ if $P \in \dyadic_{d,0}$, we get
\begin{align*}
  \big|(\ext f)_{Q,\mu_\alpha} - (\ext f)_{Q',\mu_\alpha}\big|^p & \lesssim \dashint_Q \dashint_{Q'} |\ext f(x) - \ext f(y)|^p d\mu_\alpha(x) d\mu_\alpha(y) \\
& \lesssim \Big( 2^{-k} \sum_{ \substack{ P\in \dyadic_{d,0} \\ \supp \psi_P \cap (Q\cup Q') \neq \emptyset }} \dashint_P |f| dm_d \Big)^p \\
& \lesssim 2^{-kp} \sum_{ \substack{ P\in \dyadic_{d,0} \\ \supp \psi_P \cap (Q\cup Q') \neq \emptyset }} \int_P |f|^p dm_d.
\end{align*}
Since the admissible cubes $Q$ above are relatively far away from $\real^d$, we have $\mu_\alpha(Q) \approx 2^{-k(d+1)}$, so
\[
  O^{(1)}_k \lesssim 2^{-k(d+1+p)} \sum_{Q \in \dyadic^1_k} \sum_{Q'\sim Q} \sum_{ \substack{ P\in \dyadic_{d,0} \\ \supp \psi_P \cap (Q\cup Q') \neq \emptyset }} \int_P |f|^p dm_d,
\]
and since each $P \in \dyadic_{d,0}$ appears at most some constant times $2^{(d+1)k}$ times in the above triple sum, we get
\[
  O^{(1)}_k \lesssim 2^{-kp} \sum_{P \in \dyadic_{d,0}} \int_P |f|^pdm_d = 2^{-kp} \|f\|_{L^p(\real^d,m_d)}^p.
\]
Thus,
\beqla{eq:energy-estimate-4}
  \sum_{k \geq 0} 2^{ksp} O^{(1)}_k \lesssim \sum_{k\geq 0} 2^{k(s-1)p} \|f\|_{L^p(\real^d,m_d)}^p \approx \|f\|_{L^p(\real^d,m_d)}^p.
\eeq

Now suppose that $Q \in \dyadic^2_k$, $Q' \in \dyadic^{2,*}_k$ and $Q \sim Q'$. Let $P$ and $P'$ be the (unique) cubes in $\dyadic^0_d$ such that $Q \subset \W(P)$ and  $Q' \subset \W(P')$. We evidently have $\ell(\W(P)) \geq 2^{-k}$ and $\ell(\W(P')) \approx \ell(\W(P))$. Using the Lipschitz continuity of the bump functions in the definition of $\ext f$ in conjunction with the fact that the bump functions form a partition of unity in $Q\cup Q'$, we get
\begin{align}
  \big|(\ext f)_{Q,\mu_\alpha} - (\ext f)_{Q',\mu_\alpha}\big|^p & \lesssim \dashint_{Q}\dashint_{Q'} |\ext \big( f - f_{P}\big)(x) - \ext \big( f - f_{P}\big)(y) \big|^p d\mu_{\alpha}(x)d\mu_{\alpha}(y) \notag \\
& \lesssim \frac{2^{-kp}}{\ell(P)^p} \sum_{\substack{R \in \dyadic^0_d \\ \overline{\W(R)} \cap (\overline{\W(P)}\cup  \overline{\W(P')}) \neq \emptyset} } \big|f_P - f_R \big|^p \label{eq:common-energy-estimate-3} \\
& \lesssim \frac{2^{-kp}}{\ell(P)^p}\Big( \sum_{\substack{R \in \dyadic^0_d \\ \overline{\W(R)} \cap \overline{\W(P)} \neq \emptyset} } \big|f_P - f_R \big|^p + \sum_{\substack{R \in \dyadic^0_d \\ \overline{\W(R)}\cap \overline{\W(P')} \neq \emptyset} } \big|f_{P'} - f_R \big|^p \Big).\notag
\end{align}
Since $w_\alpha \approx \ell(P)^{\alpha}$ in $\W(P)$, we have $\mu_\alpha(Q) \approx 2^{-k(d+1)}\ell(P)^{\alpha}$, so
\begin{align*}
  \mu_\alpha(Q) \big|(\ext f)_{Q,\mu_\alpha} - (\ext f)_{Q',\mu_\alpha}\big|^p & \lesssim 2^{-k(d+1+p)} \Big( \ell(P)^{\alpha-p} \sum_{\substack{R \in \dyadic^0_d \\ \overline{\W(R)} \cap \overline{\W(P)} \neq \emptyset} } \big|f_P - f_R \big|^p \\
& \qquad + \ell(P')^{\alpha-p} \sum_{\substack{R \in \dyadic^0_d \\ \overline{\W(R)} \cap \overline{\W(P')} \neq \emptyset} } \big|f_{P'} - f_R \big|^p \Big).
\end{align*}
Now summing over admissible $Q$ and $Q'$, geometric considerations imply that the terms $P\in \dyadic^0_d$ and $P' \in \dyadic^0_d$ (with $\ell(P) \geq 2^{-k}$ and $\ell(P') \geq 2^{-k}$) will appear at most a constant times $(2^{k}\ell(P))^{d+1}$ times in the resulting triple sum, so
\begin{align*}
  & \sum_{Q \in \dyadic^{2}_k}\mu_\alpha(Q) \sum_{\substack{Q' \sim Q \\ Q' \in \dyadic^{2,*}_k}} \big|(\ext f)_{Q,\mu_\alpha} - (\ext f)_{Q',\mu_\alpha}\big|^p\\
& \qquad \qquad \lesssim 2^{-kp} \sum_{\substack{P \in \dyadic^{0}_d\\ \ell(P) \geq 2^{-k} }} \ell(P)^{d+1+\alpha-p} \sum_{\substack{R \in \dyadic^0_d \\ \overline{\W(R)} \cap \overline{\W(P)} \neq \emptyset} } \big|f_P - f_R \big|^p \\
& \qquad \qquad \lesssim 2^{-kp} \sum_{0 \leq n \leq k} 2^{-n(d+1+\alpha-p)} \sum_{P \in \dyadic_{d,n}} \sum_{\substack{R \in \dyadic^0_d \\ \overline{\W(R)} \cap \overline{\W(P)} \neq \emptyset} } \big|f_P - f_R \big|^p \\
& \qquad \qquad = 2^{-kp} \sum_{0 \leq n \leq k} 2^{-n(1+\alpha-p)} \sum_{P \in \dyadic_{d,n}} m_d(P) \sum_{\substack{R \in \dyadic^0_d \\ \overline{\W(R)} \cap \overline{\W(P)} \neq \emptyset} } \big|f_P - f_R \big|^p \\
& \qquad \qquad =: 2^{-kp} \sum_{0 \leq n \leq k} 2^{-n(1+\alpha-p)} O'_n.
\end{align*}
Multiplying this by $2^{ksp}$ and summing over $k \geq 0$, we get
\beqla{eq:energy-estimate-5}
  \sum_{k \geq 0} 2^{ksp} O^{(2)}_k \lesssim \sum_{n \geq 0} 2^{-n(1+\alpha-p)} O'_n \sum_{k \geq n} 2^{k(s-1)p} \approx \sum_{n \geq 0}2^{n(s - \frac{\alpha+1}p)p}O'_n \lesssim \|f\|_{\bptd}^p,
\eeq
where the last estimate follows from the definition of the norm.

Finally, let us consider the terms in the sum $O^{(3)}_k$. Let $Q \in \dyadic^{3}_k$ and $Q' \sim Q$. Define $P := P_Q \in \dyadic_{d,k}$ as the projection of $Q$ on $\real^d$, and let $P'$ be a neighbor of $P$ in $\dyadic_{d}$ -- we will specify the choice of $P'$ later. We have
\begin{align}
&  \mu_{\alpha}(Q) \big|(\ext f)_{Q,\mu_\alpha} - (\ext f)_{Q',\mu_\alpha}\big|^p \notag\\
& \qquad \lesssim \int_{Q} \big| \ext f - f_P \big|^p d\mu_\alpha + \int_{Q'} \big| \ext f - f_{P'} \big|^p d\mu_\alpha + \mu_{\alpha}(Q) \big|f_{P} - f_{P'}\big|^p. \label{eq:energy-estimate-6}
\end{align}
To estimate the first integral above, note that
\begin{align*}
  \int_{Q} \big| \ext f - f_P \big|^p d\mu_\alpha & = \sum_{\substack{R \in \dyadic_d \\ R \subset P}} \int_{\W(R)} \big| \ext f - f_P \big|^pd\mu_\alpha \\
& = \sum_{n \geq k} \sum_{\substack{R \in \dyadic_{d,n} \\ R \subset P}} \int_{\W(R)} \big| \ext f - f_P \big|^p d\mu_\alpha.
\end{align*}
For $R \in \dyadic_{d,n}$ as in the sum above, denote by $R^{(j)}$, $k \leq j \leq n$, the (unique) cube in $\dyadic_{d,j}$ that contains $R$. Taking $\epsilon \in (0,1+\alpha)$, we get
\begin{align*}
  \big| \ext f(x) - f_P \big|^p & \lesssim \big| \ext f(x) - f_R \big|^p + \Big(\sum_{j=k+1}^n \big| f_{R^{j}} - f_{R^{j-1}} \big| \Big)^p \\
& \lesssim \sum_{j=k}^{n} 2^{(n-j)\epsilon} \sum_{\substack{R' \in \dyadic_{d,j} \\ R' \supset R}} \sum_{R'' \sim R'} \big| f_{R'} - f_{R''} \big|^p \\
& \approx \sum_{j=k}^{n} 2^{(n-j)\epsilon} 2^{jd} \sum_{\substack{R' \in \dyadic_{d,j} \\ R' \supset R}} m_d(R') \sum_{R'' \sim R'} \big| f_{R'} - f_{R''} \big|^p \\
& = \ell(R)^{-\epsilon} \sum_{\substack{R' \in \dyadic_{d} \\ R \subset R' \subset P}} \ell(R')^{-d+\epsilon} m_d(R') \sum_{R'' \sim R'} \big| f_{R'} - f_{R''} \big|^p
\end{align*}
and since $\mu_\alpha(\W(R)) \approx \ell(R)^{d+1+\alpha}$, we arrive at
\begin{align*}
  \int_{Q} \big| \ext f - f_P \big|^p d\mu_\alpha & \lesssim \sum_{\substack{R \in \dyadic_d \\ R \subset P}} \ell(R)^{d+1+\alpha-\epsilon} \sum_{\substack{R' \in \dyadic_{d} \\ R \subset R' \subset P}} \ell(R')^{-d+\epsilon} m_d(R') \sum_{R'' \sim R'} \big| f_{R'} - f_{R''} \big|^p \\
& = \sum_{\substack{R' \in \dyadic_{d} \\ R' \subset P}} \Big( \ell(R')^{-d+\epsilon} m_d(R') \sum_{R'' \sim R'} \big| f_{R'} - f_{R''} \big|^p\Big) \sum_{\substack{R \in \dyadic_d \\ R \subset R'}}  \ell(R)^{d+1+\alpha-\epsilon}.
\end{align*}
Geometric considerations again imply that every $\ell(R) \in \{ \ell(R'), \ell(R')/2, \ell(R')/4, \cdots \}$ in the innermost appears $(\ell(R')/\ell(R))^d$ times, and since $1+\alpha-\epsilon > 0$, the sum in question is comparable to $\ell(R')^{d+1+\alpha-\epsilon}$. Thus,
\beqla{eq:energy-estimate-7}
  \int_{Q} \big| \ext f - f_P \big|^p d\mu_\alpha \lesssim \sum_{\substack{R' \in \dyadic_{d} \\ R' \subset P}} \ell(R')^{1+\alpha} m_d(R') \sum_{R'' \sim R'} \big| f_{R'} - f_{R''} \big|^p.
\eeq
Now to estimate the second term in \eqref{eq:energy-estimate-6}, we have to specify the choice of $P'$. If $\overline{Q'} \cap \real^d\times\{0\} \neq \emptyset$, we define $P$ analogously to $P'$, and the integral in question can be estimated by the right-hand side of \eqref{eq:energy-estimate-7}, with $Q$ replaced $Q'$ and $P$ replaced by $P'$. If on the other hand $\overline{Q'} \cap \real^d\times\{0\} = \emptyset$, we can take $P' \in \dyadic_{d,k}\cup\dyadic_{d,k+1}$ so that $Q' = \W(P')$, which yields
\[
  \int_{Q'} \big| \ext f - f_{P'} \big|^pd\mu_\alpha \lesssim \mu_\alpha(Q') \sum_{P'' \sim P'} \big| f_{P'} - f_{P''} \big|^p \approx \ell(P')^{1+\alpha} m_d(P') \sum_{P'' \sim P'} \big| f_{P'} - f_{P''} \big|^p.
\]
Finally, the estimate for the third term in \eqref{eq:energy-estimate-6} is obvious:
\beqla{eq:energy-estimate-8}
  \mu_{\alpha}(Q) \big|f_{P} - f_{P'}\big|^p \approx \ell(P)^{1+\alpha} m_d(P)  \big|f_{P} - f_{P'}\big|^p.
\eeq
Putting together \eqref{eq:energy-estimate-7}, \eqref{eq:energy-estimate-8} and a suitable estimate for the second term in \eqref{eq:energy-estimate-6}, we get
\[
  \mu_{\alpha}(Q) \sum_{Q' \sim Q} \big|(\ext f)_{Q,\mu_\alpha} - (\ext f)_{Q',\mu_\alpha}\big|^p \lesssim \sum_{\substack{R' \in \dyadic_{d} \\ R' \subset P_Q^*}} \ell(R')^{1+\alpha} m_d(R') \sum_{R'' \sim R'} \big| f_{R'} - f_{R''} \big|^p,
\]
where $P_Q^* := P \cup \bigcup_{P' \sim P} P'$. Since each $R' \in \dyadic_d$ (with $\ell(R') \leq \min(2^{-k+1},1)$) is contained in a finite number of admissible cubes $P_Q^*$, we thus have
\[
  O^{(3)}_k \lesssim \sum_{n\geq (k-1)_+} 2^{-n(1+\alpha)} \sum_{R' \in \dyadic_{d,n}} m_d(R') \sum_{R'' \sim R'} \big| f_{R'} - f_{R''} \big|^p =: \sum_{n\geq (k-1)_+} 2^{-n(1+\alpha)} O''_n,
\]
and so
\beqla{eq:energy-estimate-9}
  \sum_{k \geq 0} 2^{ksp} O^{(3)}_k \lesssim \sum_{n \geq 0} 2^{-n(1+\alpha)} O''_n \sum_{0 \leq k \leq n+1} 2^{ksp} \approx \sum_{n \geq 0} 2^{n(s - \frac{1+\alpha}p)p} O''_n \approx \|f\|_{\bptd}^p.
\eeq
Combining the estimates \eqref{eq:lp-estimate-2}, \eqref{eq:energy-estimate-4}, \eqref{eq:energy-estimate-5} and \eqref{eq:energy-estimate-9}, we finally get
\[
  \| \ext f \|_{\bpad}^p \lesssim \| \ext f \|_{L^p(\real^{d+1}_+,\mu_\alpha)}^p + \sum_{k \geq 0} 2^{ksp}\big(O^{(1)}_k + O^{(2)}_k + O^{(3)}_k \big) \lesssim \|f\|_{\bptd}^p.
\]

\smallskip
(ii) Now let $f \in \bpad$, and for $k \in \nanu_0$ write
\[
  T_k f := \sum_{Q \in \dyadic_{d,k}} \Big(\dashint_{N(Q)} f d\mu_\alpha \Big) \chi_Q,
\]
where $N(Q) = Q\times(0,\ell(Q)] \in \dyadic_{d+1,k}$ for all $Q \in \dyadic^0_{d}$. The operators $T_k$ will play a role similar to that of the operators $\T_k$ in the proof of Theorem \ref{th:sobolev-trace}.

We first show that the limit $\lim_{k \to \infty} T_k f$ exists pointwise $m_d$-almost everywhere in $\real^d$ by estimating the $L^p(\real^d)$-norm of the function
\[
  f^* := \sum_{k \geq 0} \big| T_{k+1} f - T_k f \big| + \big|T_0 f\big|.
\]
Now if $P \in \dyadic_{d,0}$, the definition of $T_0$ shows that
\begin{align}
  \int_{P} |f^*(x)|^p dm_d & \leq \int_P \Big( \sum_{k \geq 0} \Big| T_{k+1} f(x) - T_k f(x) \Big| \Big)^p dm_d (x) + \int_{N(P)} |f|^p d \mu_\alpha \notag \\
& \lesssim \sum_{k \geq 0} 2^{k\epsilon} \int_{P} \big| T_{k+1} f(x) - T_k f(x) \big|^p dm_d(x) + \int_{N(P)} |f|^p d \mu_\alpha \label{eq:lp-estimate-1},
\end{align}
where $\epsilon := sp - \alpha - 1 > 0$. To estimate the $k$-th integral above, note that
\begin{align*}
  \int_{P} \big| T_{k+1} f(x) - T_k f(x) \big|^p dm_d(x) & = \sum_{\substack{Q \in \dyadic_{d,k}\\Q\subset P}} \int_{Q} \big| T_{k+1} f(x) - T_k f(x) \big|^p dm_d(x) \\
& \lesssim \sum_{\substack{Q \in \dyadic_{d,k}\\Q\subset P}} m_d(Q) \sum_{\substack{Q' \in \dyadic^+_{d+1}\\ Q'\sim N(Q)}} \big|f_{N(Q),\mu_\alpha} - f_{Q',\mu_\alpha} \big|^p \\
& = 2^{-kd} \sum_{\substack{Q \in \dyadic_{d,k}\\Q\subset P}} \sum_{\substack{Q' \in \dyadic^+_{d+1}\\ Q'\sim N(Q)}} \big|f_{N(Q),\mu_\alpha} - f_{Q',\mu_\alpha} \big|^p.
\end{align*}
By this and \eqref{eq:lp-estimate-1}, we can estimate $\|f^*\|_{L^p(\real^d)}^p$ by
\begin{align}
 & \sum_{k \geq 0} 2^{k(\epsilon-d)} \sum_{Q \in \dyadic_{d,k}} \sum_{\substack{Q' \in \dyadic^+_{d+1}\\ Q'\sim N(Q)}} \big|f_{N(Q),\mu_\alpha} - f_{Q',\mu_\alpha} \big|^p + \|f\|_{L^p(\real^{d+1}_+,\mu_\alpha)}^p \notag \\
\approx & \sum_{k \geq 0} 2^{k(\alpha + 1 + \epsilon)} \sum_{Q \in \dyadic_{d,k}} \mu_\alpha\big( N(Q) \big) \sum_{\substack{Q' \in \dyadic^+_{d+1}\\ Q'\sim N(Q)}} \big|f_{N(Q),\mu_\alpha} - f_{Q',\mu_\alpha} \big|^p + \|f\|_{L^p(\real^{d+1}_+,\mu_\alpha)}^p \notag \\
= & \sum_{k \geq 0} 2^{ksp} \sum_{Q \in \dyadic_{d,k}} \mu_\alpha\big( N(Q) \big) \sum_{\substack{Q' \in \dyadic^+_{d+1}\\ Q'\sim N(Q)}} \big|f_{N(Q),\mu_\alpha} - f_{Q',\mu_\alpha} \big|^p + \|f\|_{L^p(\real^{d+1}_+,\mu_\alpha)}^p \label{eq:lp-boundary-estimate} \\
\lesssim & \|f\|_{\bpad}^p. \notag
\end{align}
This shows that $f^* < \infty$ pointwise $m_d$-almost everywhere, so that the limit $\res f := \lim_{k \to \infty} T_k f$ exists at these points. We may abuse notation by writing $f$ for $\res f$ in the remainder of this proof. Since $|f| \leq f^*$ pointwise $m_d$-almost everywhere, the estimate above plainly implies
\[
  \|f\|_{L^p(\real^d)} \lesssim \|f\|_{\bpad}.
\]

Now to estimate the $\bptd$-energy of $f$, let $k \in \nanu_0$ and recall that, by a calculation similar to \eqref{eq:balls-at-boundary}, $m_d(Q) / \mu_\alpha(N(Q)) \approx 2^{k(\alpha+1)}$ for all $Q \in \dyadic_{d,k}$. The estimate \eqref{eq:common-energy-estimate-2} in the proof of Theorem \ref{th:sobolev-trace} (with $T$ and $N$ in place of $\T$ and $\N$ respectively) yields
\begin{align*}
  & \sum_{Q \in \dyadic_{d,k}} m_d(Q) \sum_{Q' \sim Q} \big|f_{Q} - f_{Q'}\big|^p \\
\lesssim & \int_{\real^d} \big| f(x) - T_k f(x) \big|^p dm_d(x) + 2^{k(\alpha+1)} \sum_{Q \in \dyadic_{d,k}} \mu_\alpha\big(N(Q)\big)\sum_{Q' \sim Q} \big|f_{N(Q),\mu_\alpha} - f_{N(Q'),\mu_\alpha}\big|^p \\
& =: I_k + 2^{k(\alpha+1)} O_k,
\end{align*}
so that
\beqla{eq:energy-estimate-1}
  \sum_{k \geq 0} 2^{k(s-\frac{\alpha+1}{p})p} \sum_{Q \in \dyadic_{d,k}} m_d(Q) \sum_{Q' \sim Q} \big|f_{Q} - f_{Q'}\big|^p \lesssim \sum_{k \geq 0} 2^{k(s-\frac{\alpha+1}{p})p} I_k + \sum_{k \geq 0} 2^{ksp} O_k.
\eeq
We have
\beqla{eq:energy-estimate-2}
  \sum_{k \geq 0} 2^{ksp} O_k \lesssim \|f\|_{\bpad}^p
\eeq
by definition. To estimate the terms $I_k$, take $\epsilon \in (0,sp-\alpha-1)$ and proceed as in the estimates following \eqref{eq:lp-estimate-1} to obtain
\begin{align*}
  I_k & \lesssim \sum_{n \geq k} 2^{(n-k)\epsilon} \int_{\real^d} \big| T_{n+1}f(x) - T_n f(x) \big|^p dm_d(x) \\
& \lesssim \sum_{n \geq k} 2^{(n-k)\epsilon}  2^{-nd} \sum_{Q \in \dyadic_{d,n}} \sum_{\substack{Q' \in \dyadic^+_{d+1}\\ Q'\sim N(Q)}} \big|f_{N(Q),\mu_\alpha} - f_{Q',\mu_\alpha} \big|^p \\
& \approx \sum_{n \geq k} 2^{(n-k)\epsilon}  2^{n(\alpha+1)} \sum_{Q \in \dyadic_{d,n}} \mu_\alpha\big( N(Q) \big)\sum_{\substack{Q' \in \dyadic^+_{d+1}\\ Q'\sim N(Q)}} \big|f_{N(Q),\mu_\alpha} - f_{Q',\mu_\alpha} \big|^p \\
 & =: \sum_{n \geq k} 2^{(n-k)\epsilon}  2^{n(\alpha+1)} O'_n,
\end{align*}
so that
\begin{align}
  \sum_{k \geq 0} 2^{k(s-\frac{\alpha+1}{p})p} I_k & \lesssim \sum_{n \geq 0} 2^{n(\alpha+1 + \epsilon)} O'_n \sum_{0 \leq k \leq n} 2^{k(sp - \alpha - 1 - \epsilon) } \approx \sum_{n \geq 0} 2^{nsp} O'_n \notag
\\ & \lesssim \|f\|_{\bpad} \label{eq:energy-estimate-3}^p,
\end{align}
where the last estimate again follows from the definition of the norm. Plugging \eqref{eq:energy-estimate-2} and \eqref{eq:energy-estimate-3} into \eqref{eq:energy-estimate-1} leads to the desired energy estimate.

\smallskip
(iii) We plainly have $\res (\ext f) = f$ for all $f \in \bptd$.

As in the proof of Theorem \ref{th:sobolev-trace}, the remaining question is whether the trace operator $\res$ constructed above is of the correct form. We again refer to Subsection \ref{ss:lebesgue} in the Appendix for details on this.
\end{proof}

\section{Proof of Theorem \ref{th:triebel-trace}}\label{se:triebel}
Let us recall that in the proof of Theorem \ref{th:besov-trace}, $N(Q)$ for $Q \in \dyadic^0_d$ was defined as $Q\times (0, \ell(Q)] \in \dyadic^0_{d+1}$. Before giving the proof of Theorem \ref{th:triebel-trace}, let us introduce the auxiliary seminorm $[f]_{s,p,q,\alpha}$, defined by
\[
  [f]_{s,p,q,\alpha}^p=\int_{\real^{d+1}_+}\Big(\sum_{k=0}^\infty2^{ksq} \sum_{P\in \dyadic_{d, k}}\sum_{\substack{Q'\in\dyadic^0_{d+1}\\ Q'\sim N(P)}} \big|f_{N(P),\mu_\alpha} - f_{Q',\mu_\alpha}\big|^q\chi_{N(P)}(x)\Big)^{p/q}\, d\mu_\alpha(x),
\]
where $f\in \locint(\real^{d+1},\mu_\alpha)$ and the parameters $p$, $q$, $s$ and $\alpha$ as in the statement of Theorem \ref{th:triebel-trace}. We obviously have $[f]_{s,p,q,\alpha}\leq\|f\|_{\F^s_{p, q}(\real^{d+1},\mu_\alpha)}$ and
\[
  [f]_{s,p,p,\alpha}^p=\sum_{k=0}^\infty2^{ksp}\sum_{P\in \dyadic_{d, k}} \mu_\alpha\big(N(P)\big)\sum_{\substack{Q'\in\dyadic^0_{d+1}\\ Q'\sim N(P)}}\big|f_{N(P),\mu_\alpha} - f_{Q',\mu_\alpha}\big|^p
\]
for all admissible values of the parameters. We shall omit $\alpha$ from the notation and write $[f]_{s,p,q}$ if there is no risk of confusion.

For the proof of Theorem \ref{th:triebel-trace}, we shall need the following lemma concerning the seminorms $[f]_{s,p,q}$.
\lem{le:energy-approx}
Suppose that $0<s<1$, $1\leq p<\infty$, $0<q,\, q'\leq\infty$ and $\alpha > -1$. Then for any $f\in \locint(\real^{d+1}_+,\mu_\alpha)$, we have
\[
  [f]_{s,p,q,\alpha} \approx [f]_{s,p,q',\alpha}
\]
with the implied constants independent of $f$.
\elem
\begin{proof}
It suffices to consider the case $q' = p$. First, in order to estimate $[f]_{s,p,q}$ from above,
write
\[
  D(P) := D(f,P) := \sum_{Q'\sim N(P)} \big|f_{N(P),\mu_\alpha} - f_{Q',\mu_\alpha}\big|
\]
for $P \in \dyadic^0_d$, so that
\[
  [f]_{s,p,q}^p \approx \int_{\real^{d+1}_+}\Big(\sum_{k=0}^\infty2^{ksq} \sum_{P\in \dyadic_{d, k}}D(P)^q\chi_{N(P)}(x)\Big)^{p/q}\, d\mu_\alpha(x)
\]
(because the sum defining $D(P)$ is uniformly finite).

Note that $\bigcup_{P\in \dyadic_d^0}N(P)=\real^d\times(0, 1] = \bigcup_{j \geq 1}\bigcup_{P\in \dyadic_{d,j}} \W(P)$. Moreover, from the definitions it is easily seen that for $R$, $P \in \dyadic^0_d$, we have $N(P) \cap \W(R) \neq \emptyset$ if and only if $R$ is a proper subset of $P$, and in this case also $\W(R) \subset N(P)$. Thus, taking $ \epsilon \in (0, 1+\alpha )$ and using H\"older's inequality (or the subadditivity of $t \mapsto t^{p/q}$ if $p \leq q$) leads to
\begin{align*}
[f]_{s,p,q}^p &= \sum_{j \geq 1} \sum_{R\in \dyadic_{d,j}}\int_{\W(R)}\Big(\sum_{k=0}^{j-1}2^{ksq}\sum_{\substack{P \in\dyadic_{d,k}\\ P \supset R}}D(P)^q\chi_{N(P)}(x)\Big)^{p/q}\, d\mu_\alpha(x)\\
&= \sum_{j \geq 1} \sum_{R\in \dyadic_{d,j}}\mu_\alpha\big( \W(R)\big)\Big(\sum_{k=0}^{j-1}2^{ksq}\sum_{\substack{P \in\dyadic_{d,k}\\ P \supset R}}D(P)^q\Big)^{p/q}\\
&\lesssim \sum_{j \geq 1} \sum_{R\in \dyadic_{d,j}}\mu_\alpha\big( \W(R)\big)\sum_{k=0}^{j-1}2^{(j-k)\epsilon}2^{ksp}\sum_{\substack{P \in\dyadic_{d,k}\\ P \supset R}}D(P)^p \\
& = \sum_{k \geq 0} 2^{k(sp - \epsilon)} \Big( \sum_{P \in\dyadic_{d,k}} D(P)^p \sum_{j > k}2^{j \epsilon} \sum_{\substack{ R \in \dyadic_{d,j} \\ R \subset P}} \mu_\alpha\big( \W(R)\big) \Big).
\end{align*}
As in the previous proofs, each term $\mu_\alpha( \W(R))$ in the innermost sum above is comparable to $2^{-j(d+1+\alpha)}$, and the sum has $2^{(j-k)d}$ such terms. This together with the choice of $\epsilon$ yields
\begin{align*}
  [f]_{s,p,q}^p & \lesssim \sum_{k \geq 0} 2^{k(sp - d - \epsilon)} \Big( \sum_{P \in\dyadic_{d,k}} D(P)^p \sum_{j > k}2^{j (\epsilon - \alpha - 1)}\Big) \\
& \approx \sum_{k \geq 0} 2^{k(sp - d - \alpha - 1)} \sum_{P \in\dyadic_{d,k}} D(P)^p \\
& \approx \sum_{k \geq 0} 2^{ksp} \sum_{P \in\dyadic_{d,k}} \mu_\alpha\big(N(P)\big)D(P)^p \\
& \approx [f]_{s,p,p}^p.
\end{align*}

For the other direction, write $W(P) := P \times (\frac12\ell(P),\ell(P)]$ for all $P \in \dyadic^0_d$. Note that $W(P) \subset N(P)$ and $\mu_\alpha(W(P)) \approx \mu_\alpha(N(P))$ for all $P$, and that the cubes $W(P)$ are pairwise disjoint. We get
\begin{align*}
  [f]_{s,p,p}^p & \approx \sum_{k \geq 0} 2^{ksp} \sum_{P \in \dyadic_{d,k}} \mu_\alpha\big( W(P) \big) D(P)^p\\
& = \sum_{k \geq 0} \sum_{P \in \dyadic_{d,k}} \int_{W(P)} \Big( 2^{ks q} D(P)^q \Big)^{p/q} d\mu_\alpha \\
& \leq \sum_{k \geq 0} \sum_{P \in \dyadic_{d,k}} \int_{W(P)} \Big( \sum_{ j \geq 0} 2^{js q} \sum_{Q \in \dyadic_{d,j}}D(Q)^q\chi_{N(Q)}(x) \Big)^{p/q} d\mu_\alpha(x) \\
& \leq \int_{\real^{d+1}_+} \Big( \sum_{ j \geq 0} 2^{js q} \sum_{Q \in \dyadic_{d,j}}D(Q)^q\chi_{N(Q)}(x) \Big)^{p/q} d\mu_\alpha(x) \\
& = [f]_{s,p,q}^p.\qedhere
\end{align*}
\end{proof}

\begin{proof}[Proof of Theorem \ref{th:triebel-trace}]
(i) Let us first establish the relevant norm inequality for the Whitney extension of a function $f \in \bptd$. By Theorem \ref{th:besov-trace} and Remark \ref{re:basic-properties}, it suffices to consider the case $q < p$. As in \eqref{eq:lp-estimate-2}, we again have
\[
  \|\ext  f\|_{L^p(\real^{d+1}_+, \mu_\alpha)}\lesssim \|f\|_{L^p(\real^d)}.
\]

Now for the $\F^{s}_{p,q}(\real^{d+1}_+,\mu_\alpha)$-energy of $\ext f$, it suffices to estimate
\beqla{eq:high-frequencies}
  \int_{\real^d} \Big(\sum_{k = k_0}^\infty 2^{ksq} \sum_{Q' \sim Q^x_k} \big|(\ext f)_{Q^x_k,\mu_\alpha} - (\ext f)_{Q',\mu_\alpha}\big|^q \Big)^{p/q}\, d\mu_\alpha(x),
\eeq
where $k_0 \geq 4$ is a fixed integer (we will specify the choice of $k_0$ later), since the corresponding integral with $\sum_{k=k_0}^{\infty}$ replaced by $\sum_{k=0}^{k_0-1}$ is easily estimated by
\[
  \|\ext f\|_{L^p(\real^{d+1}_+,\mu_\alpha)}^p \lesssim \|f\|_{L^p(\real^d)}^p.
\]

To this end, we divide the cubes in $\dyadic^+_{d+1}$ into several classes as in the proof of Theorem \ref{th:besov-trace}, but this time we need to consider four different cases. More precisely, for $k \geq k_0$ write $\dyadic_k^1$ for the dyadic cubes $Q$ in $\dyadic^+_{d+1}$ with edge length $2^{-k}$ such that $\dist(Q,\real^d\times\{0\}) > 2 - 2^{-k+2}$, $\dyadic_k^2$ for the cubes $Q$ with edge length $2^{-k}$ such that $2^{-k+1} < \dist(Q,\real^d\times\{0\}) \leq 2 - 2^{-k+2}$, $\dyadic^3_k$ for the cubes with edge length $2^{-k}$ such that $2^{-k} \leq \dist (Q,\real^d\times\{0\}) \leq 2^{-k+1}$ and $\dyadic^4_k$ for the cubes with edge length $2^{-k}$ whose closures intersect $\real^d\times\{0\}$. With these choices, the quantity \eqref{eq:high-frequencies} is comparable to
\begin{align*}
  \sum_{j=1}^4\int_{\real^{d+1}_+}\Big(\sum_{k=k_0}^\infty 2^{ksq} \sum_{Q\in \dyadic^j_k}\sum_{{ Q'\sim Q}} \big|(\ext f)_{Q,\mu_\alpha} - (\ext f)_{Q',\mu_\alpha}\big|^q\chi_{Q}(x)\Big)^{p/q}\, d\mu_\alpha(x)=: \sum_{j=1}^4O^j.
\end{align*}

The necessary estimates for the term $O^4$ are already contained in Lemma \ref{le:energy-approx} and Theorem \ref{th:besov-trace}:
\[
  O^4 = [\ext f]_{s,p,q}^p  \approx [\ext f]_{s,p,p}^p \lesssim \|\ext f\|_{\bpad}^p \lesssim \|f\|_{\bptd}^p.
\]

The term $O^3$ can be estimated in a similar manner as $O^4$, since the quantity $O^3$ is also essentially independent of the parameter $q$. This is because the cubes in $\bigcup_{k \geq k_0} \dyadic^2_k$ have bounded overlap.

In order to estimate $O^1$, let us specify $k_0$: it can be taken such that whenever $Q \in \dyadic^1_k$ with $k \geq k_0$ and $Q' \sim Q$, $\supp \psi_P \cap (Q\cup Q') \neq \emptyset$ can only hold for $P \in \dyadic_{d,0}$. Using this property together with the Lipschitz continuity of the bump functions $\psi_P$, we get
\begin{align*}
  \big|(\ext f)_{Q,\mu_\alpha} - (\ext f)_{Q',\mu_\alpha}\big|^q & \leq \Big(\dashint_Q \dashint_{Q'} |\ext f(x) - \ext f(y)| d\mu_\alpha(x) d\mu_\alpha(y)\Big)^q \\
& \lesssim \Big( 2^{-k} \sum_{ \substack{ P\in \dyadic_{d,0} \\ \supp \psi_P \cap (Q\cup Q') \neq \emptyset }} \dashint |f| dm_d \Big)^q \\
& \lesssim 2^{-kq} \Big(\sum_{ \substack{ P\in \dyadic_{d,0} \\ \supp \psi_P \cap (Q\cup Q') \neq \emptyset }} \int_P |f| dm_d\Big)^q=: 2^{-kq} J_{Q,Q'}^q.
\end{align*}
Take $\epsilon\in(0,1-s)$ and $q^* > 1$ so that $1/q^*+q/p=1$. Using the estimate above together with H\"older's inequality yields
\begin{align*}
&\Big(\sum_{k\geq k_0} 2^{ksq}\sum_{Q\in \dyadic^1_k}\sum_{{ Q'\sim Q}} \big|(\ext f)_{Q,\mu_\alpha} - (\ext f)_{Q',\mu_\alpha}\big|^q\chi_{Q}(x)\Big)^{p/q}\\
\lesssim & \Big(\sum_{k\geq k_0} 2^{k(s-1+\epsilon)p}\sum_{Q\in \dyadic^1_k}\sum_{{ Q'\sim Q}} J_{Q,Q'}^p\chi_{Q}(x)\Big)\Big(\sum_{k\geq k_0} 2^{-k\epsilon q q^*} \Big)^{p/(qq^*)}\\
\approx & \sum_{k\geq k_0} 2^{k(s-1+\epsilon)}\sum_{Q\in \dyadic^1_k}\sum_{{ Q'\sim Q}} J_{Q,Q'}^p\chi_{Q}(x).
\end{align*}
Hence we have 
\begin{align*}
O^1& \lesssim \sum_{k \geq k_0} 2^{k(s-1+\epsilon)p}\sum_{Q\in \dyadic^1_k}\mu_\alpha(Q)\sum_{Q' \sim Q}J_{Q,Q'}^p\\
& \lesssim \sum_{k \geq k_0} 2^{k(s-1-\epsilon)p-k(d+1)} \sum_{Q\in \dyadic^1_k} \sum_{Q' \sim Q} \sum_{ \substack{ P\in \dyadic_{d,0} \\ \supp \psi_P \cap (Q\cup Q') \neq \emptyset }} \int_P |f|^p dm_d,
\end{align*}
and since each $P \in \dyadic_{d,0}$ appears at most some constant times $2^{(d+1)k}$ times in the above triple sum, we arrive at
\begin{align*}
O^1 & \lesssim \sum_{k\geq k_0} 2^{k(s-1+\epsilon)}\sum_{P\in \dyadic_{d, 0}}\int_P |f|^p dm_d \approx \sum_{k\geq k_0} 2^{k(s-1+\epsilon)} \|f\|^p_{L^p(\real^d, m_d)}\approx \|f\|^p_{L^p(\real^d, m_d)}.
\end{align*}

Finally let us estimate $O^2$. Suppose that $Q \in \dyadic^2_k$ and $Q \sim Q'$. Since $\dist(Q,\real^d\times\{0\})>2^{-k+1}$, $\ell(Q')\leq 2\ell(Q)=2^{-k+1}$ and $\overline {Q'}\cap\overline Q\not=\emptyset$, we have $Q'\cap \real^d\times\{0\}=\emptyset$. 
As in the proof of Theorem \ref{th:besov-trace}, we can therefore take $P := P_Q$ and $P'$ to be the cubes in $\dyadic^0_d$ such that $Q \subset \W(P)$ and  $Q' \subset \W(P')$. Moreover, the definition of $\dyadic^2_k$ implies that $Q\cup Q' \subset \bigcup_{R \in \dyadic^0_d} \W(R)$, and the bump functions $\psi_P$ form a partition of unity of the latter set. As in \eqref{eq:common-energy-estimate-3}, we thus get
\[
  \big|(\ext f)_{Q,\mu_\alpha} - (\ext f)_{Q',\mu_\alpha}\big|^q \lesssim \frac{2^{-kq}}{\ell(P)^q} \sum_{\substack{R \in \dyadic^0_d \\ \overline{\W(R)} \cap (\overline{\W(P)}\cup  \overline{\W(P')}) \neq \emptyset} } \big|f_P - f_R \big|^q,
\]
and hence
\[
  \sum_{Q'\sim Q} \big|(\ext f)_{Q,\mu_\alpha} - (\ext f)_{Q',\mu_\alpha}\big|^q\lesssim \frac{2^{-kq}}{\ell(P_Q)^q}\Big( \sum_{\substack{R \in \dyadic^0_d\\ R\sim\sim P_Q} } \big|f_{P_Q} - f_R \big| \Big)^q,
\]
where the notation $R\sim\sim P_Q$ means that there exists $R'\in \dyadic_d^0$ such that $R\sim R'$ and $R'\sim P_Q$. The latter sum obviously has a uniformly finite number of terms $|f_{P_Q} - f_R|$.

In order to apply this estimate to $O^2$, note that by the definition of the $\dyadic^2_k$'s, we have
\[
  \bigcup_{k \geq k_0} \dyadic^2_k \subset \bigcup_{P \in \dyadic^0_d} \W(P),
\]
and that if a point $x$ belongs to one of the $\W(P)$'s above, we can have $\chi_Q(x) \neq 0$ for some $Q \in \bigcup_{k \geq k_0} \dyadic^2_k$ only if $Q \subset \W(P)$, and in this case also $\ell(Q) \leq \ell(P)$. Using these facts and H\"older's inequality (with $\epsilon \in (0,1-s)$ as in the estimate for $O^1$ above), we get
\begin{align*}
O^2 &\leq \sum_{j\geq 0} \sum_{P\in\dyadic_{d, j}}\int_{\W(P)}\Big(\sum_{k \geq j} 2^{ksq} \sum_{\substack{Q\in\dyadic_{d+1, k}\\Q\subset \W(P)}}\sum_{Q'\in\dyadic_Q^2} \big|(\ext f)_{Q,\mu_\alpha} - (\ext f)_{Q',\mu_\alpha}\big|^q\chi_{Q}(x)\Big)^{p/q}\, d\mu_\alpha(x)\\
&\lesssim \sum_{j\geq0} \sum_{P\in\dyadic_{d, j}}\int_{\W(P)}\Big(\sum_{k\geq j} 2^{k(s-1)q}\sum_{\substack{Q\in\dyadic_{d+1, k}\\Q\subset \W(P)}}\frac{1}{\ell(P)^q}\Big( \sum_{\substack{R \in \dyadic^0_d\\ R\sim\sim P} } \big|f_P - f_R \big| \Big)^q \chi_{Q}(x)\Big)^{p/q}\, d\mu_\alpha(x)\\
&= \sum_{j\geq0} 2^{jp} \sum_{P\in\dyadic_{d, j}}\int_{\W(P)}\Big(\sum_{k\geq j} 2^{k(s-1)q}\sum_{\substack{Q\in\dyadic_{d+1, k}\\Q\subset \W(P)}}\Big( \sum_{\substack{R \in \dyadic^0_d\\ R\sim\sim P} } \big|f_P - f_R \big| \Big)^q \chi_{Q}(x)\Big)^{p/q}\, d\mu_\alpha(x)\\
&\lesssim \sum_{j\geq0} 2^{j(1-\epsilon)p} \sum_{P\in\dyadic_{d, j}} \sum_{k\geq j} 2^{k(s-1+\epsilon)p}\sum_{\substack{Q\in\dyadic_{d+1, k}\\Q\subset \W(P)}} \mu_{\alpha}( Q)\sum_{\substack{R \in \dyadic^0_d\\ R\sim\sim P} } \big|f_P - f_R \big|^p \\
& = \sum_{j\geq0} 2^{j(1-\epsilon)p} \sum_{P\in\dyadic_{d, j}} \mu_\alpha\big( \W(P) \big) \Big( \sum_{\substack{R \in \dyadic^0_d\\ R\sim\sim P} } \big|f_P - f_R \big|^p \Big) \sum_{k\geq j} 2^{k(s-1+\epsilon)p} \\
& \approx \sum_{j\geq0} 2^{js p} \sum_{P\in\dyadic_{d, j}} \mu_\alpha\big( \W(P) \big) \sum_{\substack{R \in \dyadic^0_d\\ R\sim\sim P} } \big|f_P - f_R \big|^p \\
& \approx \sum_{j\geq0} 2^{j(s - \frac{\alpha+1}p) p} \sum_{P\in\dyadic_{d, j}} m_d(P)  \sum_{\substack{R \in \dyadic^0_d\\ R\sim\sim P} } \big|f_P - f_R \big|^p.
\end{align*}
Finally, since for each $P$ above we have $R\sim\sim P$ for a (uniformly) finite number of cubes $R$, the above quantity is easily estimated by $\|f\|_{\bptd}^p$.

Combining the estimates for $O^1$, $O^2$, $O^3$ and $O^4$ with the $L^p$-estimate for $\ext f$, we conclude that
\[
  \|\ext f\|_{\F^{s}_{p,q}(\real^{d+1}_+,\mu_\alpha)} \lesssim \|f\|_{\bptd}
\]

\smallskip
(ii) In order to establish the existence of the trace of a function $f \in \fpa$, we proceed as in the proof of Theorem \ref{th:besov-trace} (ii). Let $f \in \F^s_{p,q}(\real^{d+1}_+,\mu_\alpha)$, define $T_k f$ for $k \in \nanu_0$ as in that proof and put
\[
  f^* := \sum_{k \geq 0} \big| T_{k+1} f - T_k f \big| + \big|T_0 f\big|.
\]
By the estimate \eqref{eq:lp-boundary-estimate} and Lemma \ref{le:energy-approx}, we have
\[
  \|f^*\|_{L^p(\real^d)} \lesssim \|f\|_{L^p(\real^{d+1}_+,\mu_\alpha)} + [f]_{s,p,p} \approx \|f\|_{L^p(\real^{d+1}_+,\mu_\alpha)} + [f]_{s,p,q} \leq \|f\|_{\F^s_{p,q}(\real^{d+1}_+,\mu_\alpha)} < \infty,
\]
so the trace $\res f := \lim_{k \to \infty} T_k f$ is well-defined $m_d$-almost everywhere in $\real^d$. The estimates \eqref{eq:lp-boundary-estimate}, \eqref{eq:energy-estimate-1} and \eqref{eq:energy-estimate-3} then imply
\[
  \|\res f\|_{\bptd} \lesssim \|f\|_{L^p(\real^{d+1}_+,\mu_\alpha)} + [f]_{s,p,p} \approx \|f\|_{L^p(\real^{d+1}_+,\mu_\alpha)} + [f]_{s,p,q} \lesssim \|f\|_{\F^s_{p,q}(\real^{d+1}_+,\mu_\alpha)}
\]
which is the desired norm estimate.

\smallskip
(iii) That $\res (\ext f) = f$ for all $f \in \bptd$ again follows plainly from the definition of $\res$. Concerning the fact that $\res$ is actually of the form required by Definition \ref{de:whitney}, we again refer to Subsection \ref{ss:lebesgue} of the Appendix.
\end{proof}

\section{The trace of a weighted Hardy-Sobolev space}\label{se:hardy-sobolev}

In this section we present a refinement of the case $p = 1$ of Theorem \ref{th:sobolev-trace}, where $W^{1,1}(\real^{d+1}_+,\mu_\alpha)$ is replaced by a \emph{weighted Hardy-Sobolev space} on $\real^{d+1}_+$.

The real-variable Hardy spaces $H^p(\real^d)$, $0 < p \leq 1$, were defined for a general dimension $d$ and exponent $p$ in the seminal paper by Fefferman and Stein \cite{FeSt}. They have since been studied extensively, as many results of harmonic analysis that fail for $p \leq 1$ work for these spaces. We refer to \cite{St} for an extensive treatment of these spaces.

A localized version of the space $H^p$, better suited e.g.~for studying functions on domains, was introduced by Goldberg \cite{Go}. A variety of similar spaces, including spaces on domains, weighted spaces on domains and Sobolev-type spaces based on the $H^p$ norm, have since been studied e.g.~in \cite{Mi,StTo,Mi2,Mi3,WY,CCYY,KS}.

Let us now define the Hardy-Sobolev space relevant to us. Fix a function $\Phi \in C^{\infty}(\real^{d+1})$ such that $\supp \Phi \subset B(0,1)$ and $\int \Phi dm_{d+1} = 1$. Following Miyachi \cite{Mi,Mi2}, for $f \in \locint(\real^{d+1}_+,m_{d+1})$, define the radial maximal function $f^+\colon \real^{d+1}_+ \to [0,\infty]$ by
\[
  f^+(x) = \sup_{0 < t < \min(x_{d+1},1)} \big|(f * \Phi_t)(x)\big|,
\]
where $x_{d+1}$ is the $(d+1)$-th coordinate of $x$ and $\Phi_t := t^{-(d+1)}\Phi(\cdot/t)$. If $\mu$ is a Borel regular and absolutely continuous measure on $\real^{d+1}_+$, define the localized Hardy space $h^1(\real^{d+1}_+,\mu)$ as the space of locally $m_{d+1}$-integrable functions $f$ on $\real^{d+1}_+$ such that
\[
  \|f\|_{h^1(\real^{d+1}_+,\mu)} := \|f^+\|_{L^1(\real^{d+1}_+,\mu)}
\]
is finite. We clearly have $|f(x)| \leq |f^+(x)|$ for almost all $x$, so $h^1(\real^{d+1}_+,\mu) \subset L^1(\real^{d+1}_+,\mu)$ with a continuous embedding.

It follows from Miyachi's results (see also \eqref{eq:grand-maximal-estimate} below) that for the measures $\mu$ relevant to us, the space defined above is independent of $\Phi$ in the sense that two admissible choices yield the same space with equivalent norms. In fact, it will be convenient for us to choose $\Phi$ so that $\supp \Phi \subset B(0,1/8)$.

Now the Hardy-Sobolev space $h^{1,1}(\real^{d+1}_+,\mu)$ is defined as the space of functions $f \in \locint(\real^{d+1}_+,m_{d+1})$ such that the first-order distributional derivatives $\partial_j f$, $1 \leq j \leq d+1$, also belong to $\locint(\real^{d+1}_+,m_{d+1})$ and
\[
  \|f\|_{h^{1,1}(\real^{d+1}_+,\mu)} := \|f\|_{L^1(\real^{d+1}_+,\mu)} + \sum_{j=1}^{d+1} \|\partial_j f\|_{h^{1}(\real^{d+1}_+,\mu)}
\]
is finite.

The trace theorem for these spaces then reads as follows.

\thm{th:hardy-sobolev-trace}
Let $\alpha \in (-1,0)$ Then $\big(\B^{-\alpha}_{1, 1}(\real^d),h^{1,1}(\real^{d+1}_+,\mu_\alpha)\big)$ is a Whitney trace-extension pair.
\ethm

Before proving this Theorem, let us formulate a sampling lemma which is essentially folklore. For the convenience of the reader, a proof is presented in Subsection \ref{ss:sampling} of the Appendix.

\lem{le:sampling}
Suppose that $\Omega$ is an open subset of $\real^d$, that $\mu$ is a doubling measure on $\Omega$ such that every Euclidean ball (restricted to $\Omega$) has positive and finite $\mu$-measure and that $0 < \lambda < 1$. Then there is a constant $C$ depending only on the dimension $d$, the doubling constant of $\mu$ and $\lambda$ such that the following statement holds.

For every cube $Q \subset \Omega$ and $f \in L^1(Q)$, there exists a cube $\tilde Q \subset Q$ with $\ell(\tilde Q) = \lambda \ell(Q)$ such that
\[
  \dashint_{Q} \big| f - f_{ Q,\mu}\big| d\mu \leq C \dashint_{\tilde Q} \big| f - f_{\tilde Q,\mu}\big| d\mu.
\]
\elem

\begin{proof}[Proof of Theorem \ref{th:hardy-sobolev-trace}]
(i) In order to estimate the $h^{1,1}(\real^{d+1}_+,\mu_\alpha)$-norm of the Whitney extension of a function $f \in B^{-\alpha}_{1,1}(\real^d)$, we proceed as in the proof of Theorem \ref{th:sobolev-trace}. First, the $L^1$-norm of $\ext f$ can be estimated as in \eqref{eq:sobolev-lp-estimate}. In order to estimate the $h^1(\real^{d+1}_+,\mu_\alpha)$-norm of a partial derivative $\partial_j f$, write $X_1 := \cup_{Q \in \dyadic^0_d} \W(Q)$ and $X_2 := \real^{d+1}_+ \setminus X_1$.

Suppose first that $x \in X_1$, i.e.~$x \in \W(P)$ for some $Q$, and $0 < t < \min(x_{d+1},1)$. We plainly have
\[
  \big(\partial_j(\ext f)\big)*\Phi_t(x) = \Big( \partial_j \big( \ext f - f_P \big) \Big)*\Phi_t(x).
\]
Now since $x_{d+1} < 2\ell(Q)$ and we assumed the support of $\Phi$ to be contained in $B(0,1/8)$, we see that
\[
  \supp \Phi_t(x-\cdot) \subset \frac54 \W(Q),
\]
and hence $\supp \psi_{Q} \cap \supp \Phi_t(x-\cdot) \neq \emptyset$ can only hold if $Q \sim P$. Since also the $L^1$-norm of $\Phi_t(x-\cdot)$ does not depend on $t$, we get
\begin{align*}
  \big(\partial_j(\ext f)\big)^+(x) & \leq \sup_{0 < t < \min(x_{d+1},1)} \int_{\frac54\W(Q)} \sum_{ \substack{Q \sim P } } |f_Q - f_P| |\partial_j (\psi_Q)(y) ||\Phi_t(x-y)| dm_{d+1}(y) \\
& \lesssim \sum_{ Q \sim P } \frac{1}{\ell(Q)} |f_Q - f_P|.
\end{align*}
This is estimate corresponds to \eqref{eq:common-energy-estimate-1} in the proof of Theorem \ref{th:sobolev-trace}, so $\|(\partial_j(\ext f))^+\|_{L^1(X_1,\mu_\alpha)}$ can be estimated in the same way as in that proof.

Now if $x \in X_2$ and $0 < t < \min(x_{d+1},1)$, we can only have $\supp \psi_{Q} \cap \supp \Phi_t(x-\cdot) \neq \emptyset$ if $Q \in \dyadic_{0,d}$. Thus,
\begin{align*}
  \big|\big(\partial_j(\ext f)\big)*\Phi_t(x)\big| & \leq \int_{\real^{d+1}_+} \sum_{Q \in \dyadic_{d,0}} |f_Q| |\partial_j (\psi_Q)(y) | |\Phi_t(x-y)| dy m_{d+1} \\
& \lesssim \sum_{\substack{ Q \in \dyadic_{d,0}\\ \supp \psi_{Q} \cap \supp \Phi_t(x-\cdot) \neq \emptyset }} |f_Q| \leq \sum_{\substack{ Q \in \dyadic_{d,0}\\ \supp \psi_{Q} \cap B(x,1/8) \neq \emptyset }} |f_Q|
\end{align*}
Since the $\mu_\alpha$-measures of the $\frac18$-neighborhoods of the supports of $\psi_Q$ above are comparable to $1$, we get
\[
  \| (\partial_j(\ext f))^+\|_{L^1(X_2,\mu_\alpha)} \lesssim \sum_{Q \in\dyadic_{d,0}} |f_Q| \leq \|f\|_{L^1(\real^d)}.
\]
That $\res(\ext f) = f$ is then checked as in the previous proofs. This finishes the proof of part (ii).

\smallskip
(ii) Let us recall some notation from the proof of Theorem \ref{th:sobolev-trace}. For $Q \in \dyadic^0_d$, write $\W(Q) := Q \times [\ell(Q),2\ell(Q)) \in \dyadic^0_{d+1}$ and $\N(Q) := \frac54 \W(Q)$.

Now in order to verify the existence of the trace of a Hardy-Sobolev function and estimate its norm, the argument in part (ii) of the proof of Theorem \ref{th:sobolev-trace} applies here as well, as long as we can verify that every $f \in h^{1,1}(\real^{d+1}_+, \mu_\alpha)$ satisfies a suitable Poincar\'e-type inequality on cubes that are relatively far away from the boundary $\real^d$. More precisely, it suffices to show that there exists a measurable function $g \colon \real^{d+1}_+ \to [0,\infty]$ such that
\beqla{eq:hs-poincare-1}
  \dashint_{\N(Q)} \big| f - f_{\N(Q)} \big| dm_{d+1} \lesssim \ell(Q) \dashint_{\N(Q)} g dm_{d+1}
\eeq
for all $Q \in \dyadic^0_d$ (with the implied constant independent of $f$) and
\beqla{eq:hs-poincare-2}
  \|g\|_{L^1(\real^{d+1}_+),\mu} \lesssim \|f\|_{h^{1,1}(\real^{d+1}_+, \mu_\alpha)}.
\eeq

To this end, let us recall the definition of the \emph{grand maximal function} related to the space $h^1$. For $h \in \locint(\real^{d+1}_+)$ and $N \in \nanu$, define the function $\hlmax^*_N h\colon \real^{d+1}_{+} \to [0,\infty]$ by
\[
  \hlmax^*_N h (x) = \sup_{\psi \in {\mathcal F}_N(x)} \Big| \int_{\real^{d+1}_+} h(y) \psi(y) dm_{d+1}(y)\Big|,
\]
where
\begin{align*}
  \mathcal{F}_N(x) = & \Big\{\psi \in C^{\infty}(\real^{d+1}_+) \,:\, \text{there exist } y \in \real^{d+1}_+ \text{ and } r \in (0,1) \text{ such that } \\
& \qquad x \in B(y,r) \subset \real^{d+1}_+\text{, } \supp \psi \subset B(y,r) \text{ and } |\partial^{\beta} \psi| \leq r^{-(d+1)-|\beta|} \\
& \qquad \text{for all multi-indices }\beta\text{ such that }|\beta| \leq N\Big\}.
\end{align*}
We claim that
\[
  g := \sum_{j=1}^{d+1} \hlmax^*_1 (\partial_j f)
\]
satisfies \eqref{eq:hs-poincare-1} and \eqref{eq:hs-poincare-2}.

Now by \cite[Theorem 7]{KS}, there exists a constant $c$ depending only on the dimension $d$ such that
\[
  |f(x) - f(y)| \leq c |x-y| \big( g(x) + g(y) \big)
\]
for all $x$, $y \in \real^{d+1}_+$ such that $|x-y| < \min( x_{d+1}, y_{d+1}, 1)$. We can apply this estimate in a cube $\N(Q)$ as follows. Since $\dist(\N(Q),\real^d) \approx \ell(\N(Q))$, we can use Lemma \ref{le:sampling} to find a cube $\tilde Q \subset \N(Q)$ such that $\ell(\tilde Q) \approx \ell(\N(Q))$,
\[
  \dashint_{\N(Q)} \big| f - f_{\N(Q)} \big| dm_{d+1} \lesssim \dashint_{\tilde Q} \big| f - f_{\tilde Q} \big| dm_{d+1}
\]
and $|x-y| < \min(x_{d+1},y_{d+1},1)$ for all $x$, $y \in \tilde Q$. Thus,
\[
  \dashint_{\N(Q)} \big| f - f_{\N(Q)} \big| dm_{d+1} \lesssim \dashint_{\tilde Q} \dashint_{\tilde Q} |x-y|\big( g(x) + g(y) \big) dy dx \lesssim \ell(Q) \dashint_{\N(Q)} g dm_{d+1}, 
\]
which is \eqref{eq:hs-poincare-1}.

As for \eqref{eq:hs-poincare-2}, we denote by $\tilde g_j$, $1 \leq j \leq d+1$, the function on $\real^{d+1}$ that is obtained by extending $(\partial_j f)^+$ as zero on $\real^{d+1}\setminus \real^{d+1}_+$. Then by \cite[Corollary 2]{Mi}, there exists an exponent $q \in (0,1)$ and a constant $C$ independent of $f \in h^{1,1}(\real^{d+1}_+,\mu_\alpha)$ such that
\beqla{eq:grand-maximal-estimate}
  \hlmax^*_1(\partial_j f)(x) \leq C \Big( \hlmax (\tilde g_j^q)(x) \Big)^{1/q}
\eeq
for all $x \in \real^{d+1}_+$, where $\hlmax$ stands for the standard Hardy-Littlewood maximal operator on $\real^{d+1}$. Because $\alpha \in (-1,0)$, $w_\alpha$ can be extended in a natural way as an $A_{1/q}$-weight on $\real^{d+1}$, which in particular means that $\hlmax$ is bounded on $L^{1/q}(\real^{d+1},\mu_\alpha)$. Thus,
\[
  \|\hlmax^*_1(\partial_j f)\|_{L^1(\real^{d+1}_+,\mu_\alpha)} \lesssim \|\tilde g_j\|_{L^1(\real^{d+1},\mu_\alpha)} \approx \|\partial_j f\|_{h^1(\real^{d+1}_+,\mu_\alpha)},
\]
and summing up over $j$ yields \eqref{eq:hs-poincare-2}.

\smallskip
(iii) As in the previous proofs, we have $\res(\ext f) = f$ for all $f \in \B^{-\alpha}_{1, 1}(\real^d)$. The discussion concerning the form of the trace operator $\res$ is again postponed until Subsection \ref{ss:lebesgue} of the Appendix.
\end{proof}

\section{Appendix}\label{se:appendix}
In this section we present some details which were, for the sake of presentation, omitted in the previous sections.

\subsection{Coincidence of trace operators}\label{ss:lebesgue}

Recall that it was not a priori obvious that the trace operators constructed in the proofs of Theorems \ref{th:sobolev-trace}, \ref{th:besov-trace}, \ref{th:triebel-trace} and \ref{th:hardy-sobolev-trace} are of the form required by Definition \ref{de:whitney}. In this subsection we explain why this is the case.

Suppose that $f \in \bpa$ or $f \in \fpa$ with the parameters $p$, $q$ and $\alpha$ admissible for the trace theorems concerning these spaces. Then, because of \eqref{eq:balls-at-boundary} and the fact that the measure $\mu_\alpha$ is doubling on $\real^{d+1}_+$, we have that for $m_d$-almost all $x \in \real^d$, there exists a number $c \in \complex$ such that
\beqla{eq:lebesgue-2}
  \lim_{r \to 0} \dashint_{B((x,0),r)} |f(y) - c| d\mu_{\alpha} = 0.
\eeq
In fact, the set of points $x$ for which this does not hold has Hausdorff dimension at most $\max(d+1+\alpha-sp,0) < d$. This follows from a well-known covering argument and a Poincar\'e-type inequality for the function spaces in question; we refer to e.g.~\cite[Lemma 3.1 and Remark 3.2]{SS} for details. By the same argument and the Poincar\'e inequality established e.g.~in \cite[Theorem 4]{Bjo}, the same holds if $f \in W^{1,p}(\real^{d+1}_+,\mu_\alpha)$ and $s$ above is replaced by $1$. Finally, the aforementioned argument in \cite[Lemma 3.1 and Remark 3.2]{SS} also applies for functions $f \in h^{1,1}(\real^{d+1}_+,\mu_\alpha)$, since by a modification of the proof of \cite[Theorem 16]{KS}, $f$ has a local Haj\l asz gradient in $L^1(\real^{d+1}_+,\mu_\alpha)$, which yields a suitable $(1,1)$-Poincar\'e inequality for $f$.

From \eqref{eq:lebesgue-2}, it is then easy to see that the limits defining each trace operator in the above-mentioned proofs can be rewritten in the form \eqref{eq:lebesgue}.

\subsection{Equivalence of norms}\label{ss:equivalent-norms}

Here we present a direct proof of the equivalence of the (quasi-)norm \eqref{eq:besov-norm} with the standard Besov quasi-norm \eqref{standard-besov}.

\prop{pr:equivalent-quasinorms}
Let $\mu$ be a Borel regular measure on $\real^d$ such that every Euclidean ball has positive and finite measure, and such that $\mu$ is doubling with respect to the Euclidean metric.
If $0 < s < 1$, $1 \leq p < \infty$ and $0 < q \leq \infty$, then
\[
  \|f\|_{\bp}\approx \|f\|_{L^p(\real^d,\mu)} + \bigg( \int_{0}^\infty t^{-sq }\Big( \int_{\real^d} \dashint_{B(x,r)} |f(x) - f(y)|^p d\mu(y) d\mu(x) \Big)^{q/p} \frac{dt}{t} \bigg)^{1/q} 
\]
for all $f \in \locint(\real^d,\mu)$, where the implied constants are independent of $f$.
\eprop

\begin{proof}
Let us denote the standard Besov quasi-norm \eqref{standard-besov} by $\|f\|_{B^s_{p, q}(\real^d, \mu)}$. We first prove that $\|f\|_{\bp}\lesssim \|f\|_{B^s_{p, q}(\real^d, \mu)}$. To simplify the notation, write $dx$ for $d\mu(x)$ for the rest of this proof.

The doubling property of $\mu$ implies that $\mu(Q)\approx\mu(Q')$ if $Q$ and $Q'$ are cubes in $\dyadic_d$ with $Q\sim Q'$. 
Thus,
\begin{align*}
 \sum_{Q \in \dyadic_{d,k}}\mu(Q) \sum_{Q' \sim Q} \big|f_{Q,\mu} - f_{Q',\mu}\big|^p &\leq  \sum_{Q \in \dyadic_{d,k}}\mu(Q) \sum_{Q' \sim Q} \dashint_Q\dashint_{Q'} |f(x)-f(y)|^pdydx\\
 &\lesssim  \sum_{Q \in \dyadic_{d,k}}\sum_{Q' \sim Q}\frac{1}{\mu(Q)}\int_Q\int_{Q'} |f(x)-f(y)|^pdydx\\
 & \leq \sum_{Q \in \dyadic_{d,k}}\sum_{Q' \sim Q}\frac{1}{\mu(Q)}\int_Q\int_{B(x,C\cdot2^{-k})} |f(x)-f(y)|^pdydx\\
 &\lesssim \sum_{Q \in \dyadic_{d,k}} \int_Q\dashint_{B(x,C\cdot2^{-k})} |f(x)-f(y)|^pdydx\\
 &=\int_{\real^d} \dashint_{B(x,C\cdot2^{-k})} |f(x)-f(y)|^pdydx,
\end{align*}
where $C=4\sqrt d$ and the doubling property of $\mu$ was again used in the second-to-last line. This leads to
\begin{align*}
& \sum_{k \geq 0} 2^{ksq} \Big( \sum_{Q \in \dyadic_{d,k}}\mu(Q) \sum_{Q' \sim Q} \big|f_{Q,\mu} - f_{Q',\mu}\big|^p\Big)^q \notag \\
& \qquad \lesssim  \sum_{k = 0}^\infty 2^{ksq}\bigg(\int_{\real^d} \dashint_{B(x,C\cdot2^{-k})} |f(x)-f(y)|^pdydx\bigg)^{q/p}\notag\\
& \qquad \lesssim \sum_{k = 0}^\infty \int_{C\cdot2^{-k}}^{C\cdot2^{-k+1}}\bigg(\int_{\real^d} \dashint_{B(x,t)} |f(x)-f(y)|^pdydx\bigg)^{q/p}\frac{dt}{t^{1+sq}}\notag\\
&\qquad \leq \int_0^{\infty}\bigg(\int_{\real^d} \dashint_{B(x,t)} |f(x)-f(y)|^pdydx\bigg)^{q/p}\frac{dt}{t^{1+sq}},
\end{align*}
which implies that $\|f\|_{\bp}\lesssim \|f\|_{B^s_{p, q}(\real^d, \mu)}$.

In order to prove that $\|f\|_{B^s_{p, q}(\real^d, \mu)} \lesssim \|f\|_{\bp}$, we first note that a straightforward application of Fubini's theorem in conjunction with the doubling property of $\mu$ yields
\begin{align*}
& \int_1^{\infty}\Big(\int_{\real^d} \dashint_{B(x,t)} |f(x)-f(y)|^pdydx\Big)^{q/p}\frac{dt}{t^{1+sq}} \\
&\qquad \lesssim \int_1^{\infty}\Big(\int_{\real^d} |f(x)|^p dx + \int_{\real^d}|f(y)|^p \Big(\dashint_{B(y,t)} \frac{dx}{\mu(B(x,t))}\Big) dy\Big)^{q/p}\frac{dt}{t^{1+sq}} \\
&\qquad \approx \|f\|_{L^p(\real^d,\mu)}.
\end{align*}
To estimate the corresponding integral from $0$ to $1$, use the doubling property of $\mu$ to get
\begin{align*}
& \int_0^{1}\Big(\int_{\real^d} \dashint_{B(x,t)} |f(x)-f(y)|^pdydx\Big)^{q/p}\frac{dt}{t^{1+sq}} \\
&\qquad \lesssim \sum_{k\geq 0}\int_{2^{-k-1}}^{2^{-k}}\Big(\int_{\real^d}\dashint_{B(x, 2^{-k})}|f(x)-f(y)|^p\, dy\, dx\Big)^{q/p}\frac{dt}{t^{1+sq}}\\
&\qquad \lesssim\sum_{k\geq 0} 2^{ksq}\Big(\int_{\real^d}\dashint_{B(x, 2^{-k})}|f(x)-f(y)|^p\, dy\, dx\Big)^{q/p}\\
&\qquad =\sum_{k\geq 0} 2^{ksq}\bigg(\sum_{Q\in\dyadic_{d, k}}\int_Q\dashint_{B(x, 2^{-k})}|f(x)-f(y)|^p\, dy\, dx\bigg)^{q/p}.
\end{align*}
Let $Q\in \dyadic_{d, k}$ for some $k \geq 0$. For $x\in Q$ we obviously have $B(x, 2^{-k})\subset \bigcup_{Q'\sim Q}Q'$ and $\mu(B(x,2^{-k})) \approx \mu(Q)$. Thus,
\begin{align*}
& \int_Q\dashint_{B(x, 2^{-k})}|f(x)-f(y)|^p\, dy\, dx \\
\lesssim & \sum_{Q'\sim Q}\frac{1}{\mu(Q)}\int_Q\int_{Q'} |f(x)-f(y)|^p\, dydx.\\
\lesssim & \sum_{Q'\sim Q}\frac{1}{\mu(Q)}\int_Q\int_{Q'} |f(x)-f_{Q, \mu}|^p\, dydx
 + \sum_{Q'\sim Q}\frac{1}{\mu(Q)}\int_Q\int_{Q'} |f_{Q, \mu}-f_{Q', \mu}|^p\, dydx\\
&\qquad  +\sum_{Q'\sim Q}\frac{1}{\mu(Q)}\int_Q\int_{Q'} |f_{Q', \mu}-f(y)|^p\, dydx \\
=: & O_Q^1+O_Q^2+O_Q^3,
\end{align*}
so that
\begin{align}
 \|f\|_{B^s_{p, q}(\real^d, \mu)}&\lesssim \|f\|_{L^p(\real^d,\mu)}+\bigg(\sum_{k\geq 0} 2^{ksq} \Big(\sum_{Q\in\dyadic_{d, k}} (O_Q^1+O_Q^2+O_Q^3)\Big)^{q/p}\bigg)^{1/q}\notag\\
 &\lesssim \|f\|_{L^p(\real^d,\mu)}+\sum_{j=1,2,3}\bigg(\sum_{k\geq 0} 2^{ksq}\Big(\sum_{Q\in\dyadic_{d, k}} O_Q^j\big)^{q/p}\bigg)^{1/q}\notag\\
 &=: \|f\|_{L^p(\real^d,\mu)}+ H_1+H_2+H_3\label{HHH}.
\end{align}
We first estimate the quantity $H_2$. For each $Q \in \dyadic_{d,k}$ the doubling property yields 
$$O_Q^2=\sum_{Q'\sim Q}\mu(Q')|f_{Q, \mu}-f_{Q', \mu}|^p\approx \mu(Q)\sum_{Q'\sim Q}|f_{Q, \mu}-f_{Q', \mu}|^p,$$
and hence
\begin{equation}\label{H_2}
H_2\lesssim \bigg(\sum_{k\geq 0} 2^{ksq}\Big(\sum_{Q\in\dyadic_{d, k}} \mu(Q)\sum_{Q'\sim Q}|f_{Q, \mu}-f_{Q', \mu}|^p\Big)^{q/p}\bigg)^{1/q}\lesssim \|f\|_{\bp}.
\end{equation}
Next we estimate $H_1$. For any $x\in \real^d$ and $n \in \nanu_0$, define $Q_n^x$ as the (unique) cube in $\dyadic_{d, n}$ that contains $x$. By the Lebesgue differentiation theorem for doubling measures \cite[Theorem 1.8]{He}, we have $\lim_{n\rightarrow}f_{Q_n^x,\mu}=f(x)$ for $\mu$-almost every $x \in \real^d$. Hence, if $Q\in \dyadic_{d, k}$ and $x \in Q$, we have
$$|f(x)-f_{Q, \mu}|^p\leq\Big(\sum_{n=k}^{\infty} \big|f_{Q_n^x, \mu}-f_{Q_{n+1}^x, \mu}\big|\Big)^p \lesssim 2^{-k\epsilon}\sum_{n=k}^{\infty}2^{n\epsilon} \big|f_{Q_n^x, \mu}-f_{Q_{n+1}^x, \mu}\big|^p,$$
where $\epsilon > 0$ is chosen so that $\epsilon<sp/2$. Applying this estimate to $O_Q^1$ and using the fact that every cube has a (uniformly) finite number of neighbors, we get
\begin{align*}
O_Q^1&\lesssim \sum_{Q'\sim Q}\frac{\mu(Q')}{\mu(Q)}2^{-k\epsilon}\sum_{n=k}^{\infty}2^{n\epsilon}\int_Q \big|f_{Q_n^x, \mu}-f_{Q_{n+1}^x, \mu}\big|^p\, dx\\
&\lesssim 2^{-k\epsilon}\sum_{n=k}^{\infty}2^{n\epsilon}\sum_{\substack{Q''\in\dyadic_{d, n}\\ Q''\subset Q}}\int_{Q''} \big|f_{Q_n^x, \mu}-f_{Q_{n+1}^x, \mu}\big|^p\, dx\\
&\lesssim\sum_{n=k}^{\infty}2^{\epsilon(n-k)}\sum_{\substack{Q''\in\dyadic_{d, n}\\ Q''\subset Q}} \mu(Q'')\sum_{Q'''\sim Q''}\big|f_{Q'', \mu}-f_{Q''', \mu}\big|^p.
\end{align*}
In order to use this to estimate $H_1$, we consider two possible cases for the parameter $q$. First, if $0<q\leq p$, the subadditivity of the function $t \mapsto t^{q/p}$ and the fact that $s - \epsilon/p > 0$ yield
\begin{align}
H_1^q&\lesssim\sum_{k\geq 0}2^{ksq} \sum_{n=k}^{\infty}2^{\epsilon(n-k)q/p}\bigg(\sum_{Q\in\dyadic_{d, k}}\sum_{\substack{Q''\in\dyadic_{d, n}\\ Q''\subset Q}} \mu(Q'')\sum_{Q'''\sim Q''}\big|f_{Q'', \mu}-f_{Q''', \mu}\big|^p\bigg)^{q/p}\notag\\
&\leq \sum_{n=0}^{\infty}2^{\epsilon nq/p}\Big(\sum_{k=0}^{n}2^{kq(s-\epsilon/p)}\Big)\bigg(\sum_{{Q''\in\dyadic_{d, n}}} \mu(Q'')\sum_{Q'''\sim Q''}\big|f_{Q'', \mu}-f_{Q''', \mu}\big|^p\bigg)^{q/p}\notag\\
&\approx \sum_{n= 0}^\infty 2^{nsq}\bigg(\sum_{Q\in\dyadic_{d, n}} \sum_{Q'\sim Q}\mu(Q)|f_{Q, \mu}-f_{Q', \mu}|^p\bigg)^{q/p}\leq \|f\|_{\bp}^q. \notag
\end{align}
If on the other hand $p<q\leq\infty$, we may use H\"older's inequality and the fact that $s - 2\epsilon/p >0$, to obtain
\begin{align}
H_1^q&\lesssim\sum_{k\geq 0}2^{ksq}\bigg( \sum_{n=k}^{\infty}2^{-\epsilon(n-k)q/p}2^{2\epsilon(n-k)}\sum_{Q\in\dyadic_{d, k}}\sum_{\substack{Q''\in\dyadic_{d, n}\\ Q''\subset Q}} \mu(Q'')\sum_{Q'''\sim Q''}\big|f_{Q'', \mu}-f_{Q''', \mu}\big|^p\bigg)^{q/p}\notag\\
&\leq \sum_{k\geq 0}2^{ksq} \sum_{n=k}^{\infty}2^{2\epsilon(n-k)q/p}\bigg(\sum_{Q\in\dyadic_{d, k}}\sum_{\substack{Q''\in\dyadic_{d, n}\\ Q''\subset Q}} \mu(Q'')\sum_{Q'''\sim Q''}\big|f_{Q'', \mu}-f_{Q''', \mu}\big|^p\bigg)^{q/p}\notag\\
&\leq \sum_{n=0}^{\infty}2^{2\epsilon nq/p}\Big(\sum_{k=0}^{n}2^{kq(s-2\epsilon/p)}\Big)\bigg(\sum_{{Q''\in\dyadic_{d, n}}} \mu(Q'')\sum_{Q'''\sim Q''}\big|f_{Q'', \mu}-f_{Q''', \mu}\big|^p\bigg)^{q/p}\notag\\
&\approx \sum_{n\geq 0} 2^{nsq}\bigg(\sum_{Q\in\dyadic_{d, n}} \sum_{Q'\sim Q}\mu(Q)|f_{Q, \mu}-f_{Q', \mu}|^p\bigg)^{q/p}\leq \|f\|_{\bp}^q,\notag
\end{align}
which is the desired estimate for $H_1$. Finally, the terms $O^3_Q$ are essentially symmetric to with the terms $O^1_Q$, so $H_3$ can be estimated using the same argument as $H_1$. Combining these estimates with \eqref{H_2} and applying them to \eqref{HHH}, we arrive at
\[
  \|f\|_{B^s_{p, q}(\real^d, \mu)} \lesssim \|f\|_{\bp}. \qedhere
\]
\end{proof}

\subsection{Proof of Lemma \ref{le:sampling}}\label{ss:sampling} Here we present the proof of the sampling lemma that was used in the proof of Theorem \ref{th:hardy-sobolev-trace}.

\begin{proof}
Let $Q$ and $f$ be as in the statement. Let us first consider the case $\lambda = 3/4$. Let $Q_i \subset Q$, $1 \leq i \leq 2^d$, be the cubes with edge length $\frac34\ell(Q)$ that are situated at the corners of $Q$. Then $Q^* := \bigcap_{1 \leq i \leq 2^d} Q_i$ is a cube with edge length $\frac12 \ell(Q)$. By doubling, we get
\begin{align*}
  \dashint_{Q} \big| f - f_{Q,\mu} \big| d\mu & \lesssim \dashint_{Q} \big| f - f_{Q^*,\mu} \big| d\mu \approx \max_{1 \leq i \leq 2^d} \dashint_{Q_i} \big| f - f_{Q^*,\mu} \big| d\mu\\
&  \leq \max_{1 \leq i \leq 2^d} \Big( \dashint_{Q_i} \big| f - f_{Q_i,\mu} \big| d\mu + \big| f_{Q_i,\mu} - f_{Q^*,\mu} \big| \Big),
\end{align*}
and again using the doubling property of $\mu$ to estimate the latter term in the parentheses, we arrive at
\[
  \dashint_{Q} \big| f - f_{Q,\mu} \big| d\mu \leq c \max_{1 \leq i \leq 2^d} \dashint_{Q_i} \big| f - f_{Q_i,\mu} \big| d\mu,
\]
where the constant $c$ depends only on $d$ and the doubling constant of $\mu$.

Now suppose that $\lambda \in (0,1)$ as in the statement of the Lemma. Write $k_\lambda$ for the positive integer such that $(3/4)^{k_\lambda} \leq \lambda < (3/4)^{k_\lambda-1}$. Iterating the argument above $k_\lambda$ times yields a cube $Q^{k_\lambda} \subset Q$ such that $\ell(Q^{k_\lambda}) = (3/4)^{k_\lambda} \ell(Q)$ and
\[
  \dashint_{Q} \big| f - f_{Q,\mu} \big| d\mu \leq c^{k_\lambda} \dashint_{Q^{k_\lambda}} \big| f - f_{Q^{k_\lambda},\mu} \big| d\mu.
\]
Now one can simply take a cube $\tilde{Q} \subset Q$ that contains $Q^{k_\lambda}$ and has edge length $\lambda \ell(Q)$. By doubling, the integral on the right-hand side above can then be estimated by a constant times
\[
  \dashint_{\tilde Q} \big|f - f_{\tilde Q,\mu} \big| d\mu.\qedhere
\]
\end{proof}

\subsection{Extending functions from $\real^d$ to $\real^{d+n}$}\label{ss:higher-codim}
Here we present the generalizations of Theorems \ref{th:sobolev-trace} through \ref{th:triebel-trace} for Euclidean codimensions higher than $1$. The dimensions $d \in \nanu$ and $d + n$, $n \in \nanu$, will be fixed in the sequel. For convenience we also write $\real^d$ for $\real^d \times \{0\}^n \subset \real^{d+n}$ when there is no risk of confusion.

The spaces $W^{1, p}(\real^{d+n}, \mu)$, $\bps$ and $\fps$ are as in the Defintions \ref{de:sobolev} through \ref{de:triebel}. In what follows, we consider the measures $\mu_\alpha$, $\alpha > -n$, on $\real^{d+n}$, defined by
\[
  \mu_\alpha(E)=\int_E w_\alpha dm_{d+n},
\]
where $w_\alpha \in \locint(\real^{d+n})$ stands for the weight $x \mapsto \min(1,\dist(x,\real^d))^{\alpha}$.

In order to define the Whitney extension of a function on $\real^d$ to $\real^{d+n}$, we introduce some additional notation. For $Q \in \dyadic_{d,k}$, $k \in \integer$, define
\[
  \A_Q := \Big\{P\in \dyadic_{d+n, k}: P\subset \big( Q\times [-2^{-k+1}, 2^{-k+1}]^n\big) \setminus \big( Q\times (-2^{-k}, 2^{-k})^n \big) \Big\}
\]
It is then evident that $\#\A_Q = 4^n - 2^n \approx 1$, and that
\[
  \bigcup_{Q \in \dyadic_d} \A_Q
\]
is a Whitney decomposition of the the space $\real^{d+n}\setminus\real^d$ with respect to the boundary $\real^d$. We define the bump functions $\psi_P\colon\real^{d+n}\to[0,1]$ for all $P \in \bigcup_{Q \in \dyadic^0_d}\A_Q$ so that $\Lip \psi_P \lesssim 1/\ell(P)$, $\inf_{x \in P} \psi_P(x) > 0$ uniformly in $P$, $\supp \psi_P$ is contained in an $\ell(P)/4$-neighborhood of $P$ and
\[
  \sum_{Q\in \dyadic_d^0}\sum_{P\in \A_Q}\psi_P\equiv1 \ \ \ {\rm in}\ \ \ \bigcup_{Q\in \dyadic_d^0}\bigcup_{P\in \A_Q} P.
\]

\defin{de:whitney-codim}
(i) Let $f \in \locint(\real^d)$. Then the Whitney extension $\ext f \colon \real^{d+n} \to\complex$ is defined by
\[
  \ext f(x) = \sum_{Q \in \dyadic^0_d} \sum_{P\in \A_Q} \Big( \dashint_{Q} f dm_d \Big) \psi_P(x).
\]
This definition gives rise in the obvious way to the linear operator $\ext \colon \locint(\real^d) \to C^{\infty}(\real^{d+n})$.

(ii) Let $X \subset \locint(\real^d)$ be a quasinormed function space on $\real^d$, and let $Y$ be a quasinormed function space on the weighted space $(\real^{d+n},\mu)$. We say that $(X,Y)$ is a \emph{Whitney trace-extension pair} if they satisfy the conditions in Definition \ref{de:whitney} with $\real^{d+n}$ in place of $\real^{d+1}_+$ and with $\ext$ as defined above.
\edefin

We then have the following trace theorems.

\thm{th:sobolev-trace-codim}
Let $1\leq p<\infty$ and $-n<\alpha<p-n$. Then $\big(\B^{1-(\alpha+n)/p}_{p, p}(\real^d),W^{1, p}(\real^{d+n}, \mu_\alpha)\big)$ is a Whitney trace-extension pair.
\ethm

\thm{th:besov-trace-codim}
Let $0 < s < 1$, $1 \leq p \leq \infty$, $0 < q \leq \infty$ and $-n < \alpha < sp -n$. Then $\big(\bpts,\bpas\big)$ is a Whitney trace-extension pair.
\ethm

\thm{th:triebel-trace-codim}
Let $0 < s < 1$, $1 \leq p < \infty$, $0 < q \leq \infty$ and $-n < \alpha < sp -n$. Then $\big( \bptds,\fpas \big)$ is a Whitney trace-extension pair.
\ethm

These results can be proven by suitable modifications of the arguments in the proofs of Theorems \ref{th:sobolev-trace} through \ref{th:triebel-trace}. For the reader's convenience, we sketch the modified arguments below.

\begin{proof}[Proof of Theorem \ref{th:sobolev-trace-codim}]
(i) Let us estimate the weighted Sobolev norm of the Whitney extension of a function $f \in \B^{1-(\alpha+n)/p}_{p, p}(\real^d)$. First, if $P \in \A_Q$ for some $Q \in \dyadic^0_d$, it is easily seen that $\mu_\alpha(P) \approx \ell(Q)^\alpha m_{d+n}(P) \approx \ell(Q)^{d+n+\alpha}$. Since the supports of the bump functions $\psi_P$ in the definition of $\ext$ above have bounded overlap and $\#\A_Q \approx 1$ for all $Q \in \dyadic_d$, we get
\begin{align*}
 \int_{\real^{d+n}} |\ext f|^p d\mu_{\alpha} & \lesssim \sum_{Q \in \dyadic^0_d} \sum_{P\in \A_Q} \mu_{\alpha}(P) \dashint_{Q} |f|^p dm_d \approx \sum_{Q \in \dyadic^0_d} \ell(Q)^{\alpha+n} \int_{Q} |f|^p dm_d \lesssim \int_{\real^d} |f|^p dm_d.
\end{align*}

Now to estimate the weighted $L^p$-norm of $|\nabla (\ext f)|$, write $X_1 := \bigcup_{Q\in\dyadic^0_d} \cup_{P \in \A(Q)} P$ and $X_2 := \real^{d+n}\setminus X_1$. If $x \in X_1$, i.e.~$x \in \cup_{P \in \A(Q)}P$ for some $Q \in \dyadic^0_d$, we have $\sum_{Q' \in \dyadic^0_d} \sum_{P \in \A_{Q'}} \psi_P(x) = 1$, and the inner sum can only be nonzero for $Q' \sim Q$. Thus,
\[
  |\nabla(\ext f)(x)| \leq \sum_{Q' \sim Q} \sum_{P \in \A_Q'}|f_{Q} - f_{Q'}| |\Lip(\psi_P)(x)| \lesssim \sum_{Q' \sim Q} \frac{1}{\ell(Q)} |f_{Q} - f_{Q'}|.
\]
Since $\mu_{\alpha}(\cup_{P\in \A(Q)}) \approx \ell(Q)^{n+\alpha} m_d(Q)$, we arrive at
\[
  \int_{X_1} |\nabla(\ext f)|^p d\mu_\alpha \lesssim \sum_{Q \in \dyadic^0_d} \ell(Q)^{n+\alpha-p} m_d(Q) \sum_{Q' \sim Q} |f_Q - f_{Q'}|^p \lesssim \|f\|_{\B^{1-(\alpha+n)/p}_{p,p}(\real^d)}^p.
\]
If on the other hand $x \in X_2$, we can only have $\psi_P(x) \neq 0$ if $P \in \A_Q$ for some $Q \in \dyadic_{d,0}$, so estimating as in the part (i) of the proof of Theorem \ref{th:sobolev-trace}, we get
\[
  \int_{X_2} |\nabla(\ext f)(x)|d\mu_\alpha \lesssim \sum_{Q \in \dyadic_{d,0}} \int_Q |f|^p dm_d = \|f\|_{L^p(\real^d)}^p.
\]
Combining these estimates yields the desired norm inequality for the function $\ext f$.

\smallskip
(ii) Let us now show that the trace of a function $f \in W^{1,p}(\real^{d+n},\mu_\alpha)$ exists and estimate its Besov norm. To this end, write
\[
  \W(Q) := Q \times \big(\ell(Q),2\ell(Q)\big]^n \in \A_Q \quad \text{and} \quad \N(Q) := \frac54 \W(Q)
\]
for all $Q \in \dyadic_d$, and for $k \in \nanu_0$ write
\[
  \T_k f := \sum_{Q \in \dyadic_{d,k}} \Big( \dashint_{\N(Q)} f dm_{d+n}\Big)\chi_Q.
\]
To establish the existence of the trace function, we thus want to estimate the $L^p$-norm of the function
\[
  f^* := \sum_{k \geq 0} \big|\T_k f - \T_{k+1} f\big| + \big|\T_0 f\big|.
\]
Then, since $\mu_\alpha(\N(Q)) \approx \ell(Q)^\alpha m_{d+n}(\N(Q)) \approx \ell(Q)^{d+n+\alpha}$ for all $Q \in \dyadic^0_d$, an estimate similar to the one in the part (ii) of the proof of Theorem \ref{th:sobolev-trace} yields
\[
  \big|\T_k f(x) - \T_{k+1} f(x) \big| \lesssim 2^{-k}\Big( \dashint_{\N(Q^x_k)}|\nabla f|^pd\mu_\alpha \Big)^{1/p} + 2^{-k}\Big( \dashint_{\N(Q^x_{k+1})}|\nabla f|^pd\mu_\alpha \Big)^{1/p},
\]
and since $p - (n+\alpha) > 0$, an estimate similar to the one in the part (ii) of the proof of Theorem \ref{th:sobolev-trace} again yields
\[
  \int_{\real^d} |f^*|^p dm_d \lesssim \sum_{Q \in \dyadic^0_d} \int_{\N(Q)} |\nabla f|^p d\mu_\alpha + \sum_{P \in \dyadic_{d,0}} \int_{\N(P)} |f|^p d\mu_\alpha \lesssim \|f\|_{W^{1,p}(\real^{d+n},\mu_\alpha)}^p.
\]
Hence the trace function $\res f \in L^p(\real^d)$ exists in a suitable sense and has the correct bound for its $L^p$ norm. In the sequel, we shall simply write $f$ for $\res f$.

Now to estimate the $\B^{1-(\alpha+n)/p}_{p,p}$-energy of $f$, recall that $m_d(Q)/\mu_\alpha( \N(Q)) \approx \ell(Q)^{-(\alpha+n)}$. Hence, replacing $\alpha+1$ by $\alpha + n$ in \eqref{eq:sobolev-trace-estimate-1}, we get
\begin{align*}
  & \sum_{k \geq 0} 2^{k(1-\frac{\alpha+n}{p})p} \sum_{Q \in \dyadic_{d,k}} m_d(Q)  \sum_{Q' \sim Q} \big|f_{Q} - f_{Q'}\big|^p \\
\lesssim & \sum_{k \geq 0} 2^{k(1-\frac{\alpha+n}{p})p} \int_{\real^d} \big| f(x) - \T_k f(x) \big|^p dm_d(x) + \sum_{k \geq 0} \int_{\cup_{2^{-k-1} \leq \ell(Q') \leq 2^{-k+1}} \N(Q')} |\nabla f|^pd\mu_\alpha \\
\lesssim & \|f\|_{W^{1,p}(\real^{d+n},\mu_\alpha)}^p,
\end{align*}
which is the desired estimate.

\smallskip
(iii) As in the proofs of Theorems \ref{th:sobolev-trace} through \ref{th:triebel-trace}, it remains to verify that the trace operator $\res f$ above coincides with the one in Definition \ref{de:whitney-codim}. This again follows from the discussion in Subsection \ref{ss:lebesgue}.
\end{proof}

\begin{proof}[Proof of Theorem \ref{th:besov-trace-codim} (sketch)]
Again, we only consider the case $q = p < \infty$. In the following proof, we shall use the notation
\[
  \U_Q := \bigcup_{R \in \A_Q} R \subset \real^{d+n}
\]
for all $Q \in \dyadic_d$.

\smallskip
(i) We first establish the desired norm inequality for the extension of a function $f \in \bpts$. As in the proof of Theorem \ref{th:sobolev-trace-codim} above, we have
\[
  \|\ext f\|_{L^p(\real^{d+n})} \lesssim \|f\|_{L^p(\real^d)}.
\]
To estimate the $\bpts$-energy of $\ext f$, we divide the cubes in $\dyadic_{d+n}$ into three separate classes according to their distances to $\real^d$. For $Q \in \dyadic_{d+n}$, define
\[
  \sdist(Q,\real^d) := \inf \big\{ \max_{1 \leq i \leq d+n} |x_i - y_i| \,:\, x \in Q,\, y\in \real^d \times \{0\}^n \big\},
\]
where $x_i$ and $y_i$ stand for the $i$th coordinates of $x$ and $y$ respectively. For $k \geq 0$, write $\dyadic^1_k$ for the collection of dyadic cubes in $\dyadic_{d+n,k}$ such that $\sdist(Q,\real^d) \geq 2$, $\dyadic^2_k$ for the collection of dyadic cubes such that $2^{-k} \leq \sdist(Q,\real^d) < 2$ and $\dyadic^3_k$ for the collection of dyadic cubes whose closures intersect $\real^d$. Also write $\dyadic^{2,*}_k$ for the collectino of cubes in $\bigcup_{i = \max(k-1,0)}^{k+1}\dyadic^2_i$ that are contained in $\bigcup_{Q \in \dyadic^2_k} Q$. With these definitions, it suffices to estimate the quantity in \eqref{eq:cubes-distances} at each level $k \geq 0$.

We then have $O^{(1)}_k \lesssim \|f\|_{L^p(\real^d)}^p$ for $k \in \{0,1\}$, and for $k \geq 2$ we may estimate $O^{(1)}_k$ essentially as in part (i) of the proof of Theorem \ref{th:besov-trace}. One gets
\[
  \big|(\ext f)_{Q,\mu_\alpha} - (\ext f)_{Q',\mu_\alpha} \big|^p \lesssim 2^{-kp} \sum_{P \in \dyadic_{d,0}} \sum_{\substack{P'\in \A_{P} \\ \supp \psi_{P'} \cap (Q \cup Q') \neq \emptyset}} \int_P |f|^p dm_d
\]
for all cubes $Q \in \dyadic^1_k$ and $Q' \sim Q$. Now $\mu_\alpha(Q) \approx 2^{-k(d+n)}$ and summing the previous estimate over $Q$, each term $P'$ will appear in the resulting triple sum at most a constant times $2^{(d+n)k}$ times, so
\[
  \sum_{k \geq 0} 2^{ksp} O^{(1)}_k \lesssim \sum_{k \geq 0} 2^{k(s-1)p} \|f\|_{L^p(\real^d)}^p \approx \|f\|_{L^p(\real^d)}^p.
\]

Now to estimate the terms $O^{(2)}_k$, suppose that $Q \in \dyadic^2_k$ and $Q' \in \dyadic^{2,*}_k$ for some $k \geq 0$ and that $Q' \sim Q$. Denoting by $P$ and $P'$ the unique cubes in $\dyadic^0_d$ such that $Q \in \U_P$ and $Q' \in \U_{P'}$, the argument used in \eqref{eq:common-energy-estimate-3} yields.
\[
  \big|(\ext f)_{Q,\mu_\alpha} - (\ext f)_{Q',\mu_\alpha} \big|^p \lesssim \frac{2^{-kp}}{ \ell(P)^p } \Big(\sum_{\substack{R \in \dyadic^0_d \\ \overline{\U_R} \cap \overline{\U_P} \neq \emptyset}} \big|f_P - f_R\big|^p + \sum_{\substack{R \in \dyadic^0_d \\ \overline{\U_R} \cap \overline{\U_{P'}} \neq \emptyset}} \big|f_{P'} - f_R\big|^p \Big) 
\]
Now multiplying this estimate by $\mu_\alpha(Q) \approx 2^{-k(d+n)}\ell(P)^\alpha$ and summing over admissible $Q$ and $Q'$, it can be seen that the terms $P$ and $P'$ will appear in the resulting sum at most a constant times $(2^k \ell(P))^{d+n}$ times. Thus, the estimates for the terms $O^{(2)}_k$ in the proof of Theorem \ref{th:besov-trace} apply here as well, with $\alpha+1$ replaced by $\alpha+n$.

Finally, let $Q\in \dyadic^3_k$ and $Q' \sim Q$. Write $P := P_Q$ for the projection of $Q$ on $\real^d$, and let $P'$ be a neighbor of $P$ (to be specified later). We have
\begin{align}
&  \mu_\alpha(Q) \big| (\ext f)_{Q,\mu_\alpha} - (\ext f)_{Q',\mu_\alpha}\big|^p \notag \\
& \qquad \lesssim \int_Q |\ext f - f_P|^p d\mu_\alpha + \int_{Q'} |\ext f - f_{P'}|^p d\mu_\alpha + \mu_\alpha(Q)\big|f_P - f_{P'}\big|^p. \label{eq:o3-codim}
\end{align}
The first integral can be written as
\[
  \sum_{n \geq k} \sum_{\substack{R \in \dyadic_{d,n} \\ R \subset P}} \sum_{\substack{Q^* \in \A_R \\ Q^* \subset Q}} \int_{Q^*} \big|\ext f - f_P\big|^p d\mu_\alpha,
\]
and this sum can be estimated like the corresponding sum in the proof of Theorem \ref{th:besov-trace}, again with $1+\alpha$ replaced by $n + \alpha$, so
\[
  \int_Q |\ext f - f_P|^p d\mu_\alpha \lesssim \sum_{\substack{R' \in \dyadic_d \\ R' \subset P }} \ell(R')^{n+\alpha} m_d(R') \sum_{R'' \sim R'} \big|f_{R'} - f_{R''}\big|^p.
\]
The second term in \eqref{eq:o3-codim} can (with an appropriate choice of $P'$) be estimated either like the first term, or by
\[
  \ell(P')^{n+\alpha} m_d(P') \sum_{P'' \sim P'} |f_{P'} - f_{P''}|^p.
\]
Putting together these estimates and recalling (for the third term in \eqref{eq:o3-codim}) that $\mu_{\alpha}(Q) \approx \ell(P)^{n+\alpha} m_d(P)$, we get
\[
  \mu_\alpha(Q) \sum_{Q' \sim Q}\big| (\ext f)_{Q,\mu_\alpha} - (\ext f)_{Q',\mu_\alpha}\big|^p \lesssim \sum_{\substack{R' \in \dyadic_d \\ R' \subset P^*_Q}} \ell(R')^{n\alpha}m_d(R') \sum_{R'' \sim R'} \big|f_{R'} - f_{R''}\big|^p,
\]
where $P^*_Q := P \cup \bigcup_{P' \sim P} P'$. Part (i) of the proof can then be finished as in the proof of Theorem \ref{th:besov-trace}.

\smallskip
(ii) Now for $f \in \B^s_{p,p}(\real^{d+n}_+,\mu_\alpha)$ and $k \in \nanu_0$, write
\[
  T_k := \sum_{Q \in \dyadic_d,k} \Big( \dashint_{N(Q)} f d\mu_\alpha\Big) \chi_Q,
\]
where $N(Q) := Q \times (0,\ell(Q)]^n$, and
\[
  f^* := \sum_{k \geq 0} \big|T_{k+1}f - T_k f \big| + \big| T_0 f\big|.
\]
Repeating the corresponding argument in the proof of Theorem \ref{th:besov-trace} (with $\epsilon = sp - \alpha - n$ instead of $\epsilon = sp - \alpha - 1$), we get
\begin{align}
& \quad  \|f^*\|_{L^p(\real^d)}^p \notag\\
& \lesssim \sum_{k \geq 0} 2^{ksp} \sum_{Q \in \dyadic_{d,k}} \mu_{\alpha}\big(N(Q)\big) \sum_{\substack{Q' \in \dyadic_{d+n} \\ Q' \sim N(Q)}} \big|f_{N(Q),\mu_\alpha} - f_{Q',\mu_\alpha} \big|^p + \|f\|^p_{L^p(\real^{d+n},\mu_\alpha)} \\
& \lesssim \|f\|^p_{\real^{d+n}_+,\mu_\alpha}\notag,
\end{align}
so the trace $\res f := \lim_{k \to \infty} T_k f$ exists in $L^p(\real^d)$ and pointwise $m_d$-almost everywhere, with the correct bound for its $L^p$-norm. For the energy estimate, recall that $m_d(Q)/\mu_\alpha(N(Q)) \approx \ell(Q)^{-(\alpha+n)}$, and proceed as in the proof of Theorem \ref{th:besov-trace} (with $1+\alpha$ replaced by $n+\alpha$).

\smallskip
(iii) To see that the trace operator $\res$ constructed above can be written in the form required by Definition \ref{de:whitney-codim}, we again refer to Subsection \ref{ss:lebesgue}.
\end{proof}

For the proof of Theorem \ref{th:triebel-trace-codim}, let us introduce the sets
\[
  \N_P := \big\{ Q \in \dyadic_{d+n,k} \,\:\, \overline{Q}\cap P \neq \emptyset \big\}
\]
for all $P \in \dyadic_{d,k}$, $k \in \nanu_0$, and the quantities
\[
  \langle f \rangle_{s,p,q}^p := \langle f \rangle_{s,p,q,\alpha}^p := \int_{\real^{d+n}} \Big( \sum_{k \geq 0} 2^{ksq} \sum_{P \in \dyadic_{d,k}} \sum_{Q \in \N_P} \sum_{Q' \sim Q} \big| f_{Q,\mu_\alpha} - f_{Q',\mu_\alpha}\big|^q \chi_Q(x) \Big)^{p/q} d\mu_\alpha(x)
\]
for all $f \in \locint(\real^{d+n},\mu_\alpha)$. We then have $\langle f \rangle_{s,p,q,\alpha} \leq \|f\|_{\F^s_{p,q}(\real^{d+n}),\mu_\alpha}$ and
\[
  \langle f \rangle_{s,p,p,\alpha}^p = \sum_{k=0}^\infty \sum_{P \in \dyadic_{d,k}} \sum_{Q \in \N_P} \mu_\alpha(Q) \sum_{Q' \sim Q} \big| f_{Q,\mu_\alpha} - f_{Q',\mu_\alpha}\big|^p
\]
for all admissible values of the parameters. We also have
\beqla{eq:energy-approx-codim}
  \langle f \rangle_{s,p,q,\alpha} \approx \langle f \rangle_{s,p,q',\alpha}
\eeq
for all admissible values of the parameters, with the implied constants independent of $f$, which can be proven like Lemma \ref{le:energy-approx}.

\begin{proof}[Proof of Theorem \ref{th:triebel-trace-codim} (sketch)]
(i) In order to estimate the Triebel-Lizorkin norm of the extension of a function $f \in \bptds$, recall first that
\[
  \|\ext f\|_{L^p(\real^{d+n}_+,\mu_\alpha)} \lesssim \|f\|_{L^p(\real^d)}.
\]
For the energy estimate, it suffices to consider the quantity
\beqla{eq:high-frequencies-codim}
  \int_{\real^d} \Big( \sum_{k=k_0}^{\infty}2^{ksq} \sum_{Q \in \dyadic_{d+n,k}}\sum_{Q' \sim Q}\big|(\ext f)_{Q,\mu_\alpha}-(\ext f)_{Q',\mu_\alpha} \big|^q\chi_Q(x) \Big)^{p/q} d\mu_\alpha(x)
\eeq
with a suitably chosen $k_0 \in \nanu$ (independent of $f$). To this end, recall that the distance $\sdist(Q,\real^d)$ for $Q \in \dyadic_{d,n}$ was defined in the proof of Theorem \ref{th:besov-trace-codim} above. Now for $k \geq k_0$, write $\dyadic^1_k$ for the collection of cubes $Q$ in $\dyadic_{d+n,k}$ such that $\sdist(Q,\real^d) > 2 - 2^{-k+2}$, $\dyadic^2_k$ for the collection of cubes $Q$ in $\dyadic_{d+n,k}$ with $2^{-k+1} < \sdist(Q,\real^d) \leq 2^{-k+2}$, $\dyadic^3_k$ for the collection of cubes $Q$ in $\dyadic_{d+n,k}$ with $2^{-k} \leq \sdist(Q,\real^d) \leq 2^{-k+1}$ and $\dyadic^4_k$ for the collection of dyadic $Q$ in $\dyadic_{d+n,k}$ such that $\overline{Q}\cap\real^d \neq \emptyset$. Then \eqref{eq:high-frequencies-codim} can be estimated from above by $O^1 + O^2 + O^3 + O^4$, where each $O^j$ is defined as the quantity \eqref{eq:high-frequencies-codim} with $\dyadic^j_k$ in place of $\dyadic_{d+n,k}$ in the middle sum. As in the proof of \ref{th:triebel-trace}, it turns out that by \eqref{eq:energy-approx-codim}, the quantities $O^4$ and $O^3$ are essentially independent of the parameter $q$, so the desired norm estimate for them follows from Theorem \ref{th:besov-trace-codim}. The quantities $O^1$ and $O^2$ can be estimated by a suitable modification of the argument in the proof of Theorem \ref{th:triebel-trace}, the details being omitted.

\smallskip
(ii) To obtain the existence and norm inequality for the trace function of $f \in \F^s_{p,q}(\real^{d+n},\mu_\alpha)$, one defines $\res := \lim_{k \to \infty} T_k f$, where $T_k f$ is as in the proof of Theorem \ref{th:besov-trace-codim}, and the limit exists in $L^p(\real^d)$ with the correct norm bound. From the proof of Theorem \ref{th:besov-trace-codim} and \eqref{eq:energy-approx-codim}, one further deduces that
\[
  \|\res f\|_{\bptds} \lesssim \|f\|_{L^p(\real^{d+n},\mu_\alpha)} + \langle f \rangle_{s,p,p} \approx \|f\|_{L^p(\real^{d+n},\mu_\alpha)} + \langle f \rangle_{s,p,q} \lesssim \|f\|_{\F^s_{p,q}(\real^{d+n},\mu_\alpha)}.
\]

\smallskip
(iii) To see that the trace operator $\res$ constructed above can be written in the form required by Definition \ref{de:whitney-codim}, we again refer to Subsection \ref{ss:lebesgue}.
\end{proof}

\end{document}